\documentclass[english]{article}
\usepackage[T1]{fontenc}
\usepackage{babel}
\usepackage{refstyle}
\usepackage{units}
\usepackage{mathtools}
\usepackage{bm}
\usepackage{amsmath}
\usepackage{amsthm}
\usepackage{amssymb}
\usepackage{graphicx}
\usepackage{setspace}
\usepackage[all]{xy}
\date{}
\usepackage[pdfusetitle,
 bookmarks=true,bookmarksnumbered=false,bookmarksopen=false,
 breaklinks=false,pdfborder={0 0 1},backref=false,colorlinks=false]
 {hyperref}

\makeatletter

\AtBeginDocument{\providecommand\Lemref[1]{\ref{Lem:#1}}}
\AtBeginDocument{\providecommand\Thmref[1]{\ref{Thm:#1}}}
\AtBeginDocument{\providecommand\defref[1]{\ref{def:#1}}}
\AtBeginDocument{\providecommand\Secref[1]{\ref{Sec:#1}}}
\AtBeginDocument{\providecommand\subsecref[1]{\ref{subsec:#1}}}
\AtBeginDocument{\providecommand\thmref[1]{\ref{thm:#1}}}
\AtBeginDocument{\providecommand\Propref[1]{\ref{Prop:#1}}}
\AtBeginDocument{\providecommand\enuref[1]{\ref{enu:#1}}}
\AtBeginDocument{\providecommand\lemref[1]{\ref{lem:#1}}}
\AtBeginDocument{\providecommand\propref[1]{\ref{prop:#1}}}
\AtBeginDocument{\providecommand\figref[1]{\ref{fig:#1}}}
\AtBeginDocument{\providecommand\Subsecref[1]{\ref{Subsec:#1}}}
\AtBeginDocument{\providecommand\Figref[1]{\ref{Fig:#1}}}
\RS@ifundefined{subsecref}
  {\newref{subsec}{name = \RSsectxt}}
  {}
\RS@ifundefined{thmref}
  {\def\RSthmtxt{theorem~}\newref{thm}{name = \RSthmtxt}}
  {}
\RS@ifundefined{lemref}
  {\def\RSlemtxt{lemma~}\newref{lem}{name = \RSlemtxt}}
  {}

\numberwithin{equation}{section}
\numberwithin{figure}{section}
\usepackage{thmtools}
\newcommand{\restate}[1]{\csname #1\endcsname*}
\theoremstyle{plain}
\newtheorem{thm}{\protect\theoremname}[section]
\theoremstyle{definition}
\newtheorem{defn}[thm]{\protect\definitionname}
\newenvironment{restatableThm}[2][]
 		{\restatable[#1]{thm}{#2}}
 		{\endrestatable}
\theoremstyle{remark}
\newtheorem{rem}[thm]{\protect\remarkname}
\theoremstyle{plain}
\newtheorem{lem}[thm]{\protect\lemmaname}
\theoremstyle{remark}
\newtheorem*{rem*}{\protect\remarkname}
\theoremstyle{plain}
\newtheorem{prop}[thm]{\protect\propositionname}
\theoremstyle{plain}
\newtheorem{cor}[thm]{\protect\corollaryname}
\theoremstyle{plain}
\newtheorem{fact}[thm]{\protect\factname}

\newref{lem}{name=lemma~}
\RS@ifundefined{Propref}{\newref{prop}{name=proposition~,Name=Proposition~,names=propositions~,Names=Propositions~}}{}

\newref{sec}{name=section~,Name=Section~,names=sections~,Names=Sections~}
\newref{subsec}{name=subsection~,Name=Subsection~,names=subsection~,Names=Subsection~}

\makeatother

\usepackage{cite}

\providecommand{\corollaryname}{Corollary}
\providecommand{\definitionname}{Definition}
\providecommand{\factname}{Fact}
\providecommand{\lemmaname}{Lemma}
\providecommand{\propositionname}{Proposition}
\providecommand{\remarkname}{Remark}
\providecommand{\theoremname}{Theorem}

\begin{document}
\global\long\def\im{\text{Im}}%
\global\long\def\conj{\text{Conj}}%
\global\long\def\res{\text{Res}}%
\global\long\def\acts{\curvearrowright}%
\global\long\def\jsj{\text{JSJ}}%
\global\long\def\surject{\twoheadrightarrow}%
\global\long\def\comp{\text{Comp}}%
\global\long\def\sker{\underrightarrow{\ker}^{\omega}}%
\global\long\def\amal#1#2#3{#1\underset{#3}{*}#2}%
\global\long\def\hom{\text{Hom}}%
\global\long\def\isom{\text{Isom}}%
\global\long\def\bndv{\text{qTriv}}%
\global\long\def\bnd{\bndv_{c}}%
\global\long\def\wsurelu{\ \omega\text{-almost surely}}%
\global\long\def\onto{\twoheadrightarrow}%
\global\long\def\restriction#1#2{\left.\kern-\nulldelimiterspace#1\vphantom{\big|}\right|_{#2}}%
\global\long\def\mod{\text{Mod}}%
\global\long\def\aut{\text{Aut}}%
\global\long\def\stab{\text{Stab}}%
\global\long\def\by{\mathbf{Y}}%
 
\global\long\def\projcopx{\mathcal{P}_{K}}%
\global\long\def\diam{\text{diam}}%
\global\long\def\qtms{\mathcal{C}_{K}}%
\global\long\def\c{\mathcal{C}}%

\title{On Equational Noetherianity of Colorable Hierarchically Hyperbolic
Groups}
\author{Barak Ohana}
\maketitle
\begin{abstract}
We study the homomorphisms from a fixed finitely generated group to
strictly acylindrical colorable hierarchically hyperbolic groups.
We prove that any such group is equationally noetherian.
\end{abstract}
\tableofcontents{}

\newpage{}

\section{Introduction}

\global\long\def\hom{\text{Hom}}%

In this paper we are concerned with understanding solution sets of
systems of equations over a group $\Gamma$, i.e. algebraic sets over
$\Gamma$. Understanding those sets is the first step of understanding
the elementary theory of $\Gamma$ \cite{Se1,Se2,BPP}.

An \textit{equation} is an element $W\in\mathbb{F}_{n}=\mathbb{F}\left<x_{1},\cdots,x_{n}\right>$,
the free group on $n$ elements. For $\left(g_{1},\ldots,g_{n}\right)\in\Gamma^{n}$
we define $W\left(g_{1},\ldots,g_{n}\right)\in\Gamma$ by substituting
$g_{i}$ for $x_{i}$ (i.e. as the image of $W$ by the map $\mathbb{F}_{n}\to\Gamma$
defined by $x_{i}\mapsto g_{i}$). The tuple $\left(g_{1},\ldots,g_{n}\right)$
is a \textit{solution to $W=1$ (or just a solution to $W$, for short)
if} $W\left(g_{1},\ldots,g_{n}\right)=1$. Finally, for a system of
equations $\Sigma\subset\mathbb{F}_{n}$ we define the solution set
\[
V_{\Gamma}\left(\Sigma\right)=\left\{ \left(g_{1},\ldots g_{n}\right)\in\Gamma^{n}\mid W\left(g_{1},\ldots g_{n}\right)=_{\Gamma}1\ \forall W\in\Sigma\right\} .
\]

There is an obvious bijection between the tuples in $V_{\Gamma}\left(\Sigma\right)$
and $\hom\left(\nicefrac{\mathbb{F}_{n}}{\left<\left<\Sigma\right>\right>},\Gamma\right)$,
and it is more convenient to study the latter. Since any finitely
generated group can be represented as a quotient of the free group,
this analysis is equivalent to studying $\hom\left(G,\Gamma\right)$
for an arbitrary finitely generated group $G$.

In classical algebraic geometry, given a set of polynomial equations
$S$ over some field $k$, there exists a finite subset $S_{0}$ of
$S$ such that solving $S_{0}$ is equivalent to solving $S$. This
property, known as noetherianity, is the essence of Hilbert's basis
theorem. As a group-theoretic counterpart, we have the notion of equationally
noetherian groups.
\begin{defn}
A group $\Gamma$ is called \textbf{equationally noetherian }if for
any $n$ and any system of equations $\Sigma\subseteq\mathbb{F}_{n}$
there exists a finite subset $\Sigma_{0}\subseteq\Sigma$ such that
\[
V_{\Gamma}\left(\Sigma\right)=V_{\Gamma}\left(\Sigma_{0}\right)
\]
\end{defn}

Equational noetherianity plays a crucial role in the examination of
equations over groups and serves as a fundamental property in the
study of the first-order theory of groups (see \cite{Se1,Se7,JA-SE,Se10}).
Examples of equationally noetherian groups include $\left(i\right)$
Finite groups; $\left(ii\right)$ Linear groups; and $\left(iii\right)$
Hyperbolic groups. Moreover, equationally noetherian groups are closed
under direct and free products \cite{Se10,GR-HU}, as well as under
quotienting by algebraic normal groups and taking finite index extensions
\cite{BMR}. In \Lemref{en_is_eq_factoring} we state an equivalent
formulation of equational noetherianity using the language of homomorphisms.

\subsubsection*{Hierarchically hyperbolic groups}

Hierarchically hyperbolic groups (or HHGs for short) are a generalization
of hyperbolicity that allows ``product regions'', which are quasi-embedded
copies of finite products of simpler (i.e. lower dimension/complexity)
HHGs, where HHGs of complexity $1$ are hyperbolic spaces. The motivation
behind defining HHGs is the attempt of axiomatizing the similarities
observed between the connection of mapping class group and the curve
complex (\cite{MM00,MM99}) and the connection of cubical groups and
hyperplanes. The core idea of HHGs is the existence of hierarchy of
hyperbolic spaces encoding the HHG's entire geometry, thus allowing
to study HHGs through their relation to this hierarchy of hyperbolic
space. For a comprehensive yet introductory survey of the subject,
we direct the reader to \cite{S}. For a detailed definition, we recommend
consulting \cite{BHS1,BHS2}.

Examples of HHGs include:
\begin{itemize}
\item Hyperbolic groups and direct product of hyperbolic groups;
\item Mapping class group or finite type hyperbolic surfaces \cite{MM00,MM99,Beh06,BKMM12};
\item Certain CAT(0) cubulated groups \cite{Ha20};
\item The Teichmüller space $\mathcal{T}\left(S\right)$ of finite type
surfaces with the Weil-Peterson metric \cite{MM99,MM00,Bro03,Beh06,BKMM12};
\item The Teichmüller space $\mathcal{T}\left(S\right)$ with the Teichmüller
metric \cite{MM99,Ra07,Dur,EMR};
\end{itemize}
In this paper we show that all HHGs that satisfy certain properties
are equationally noetherian.
\begin{restatableThm}{main-thm}
\label{thm:main-thm}Every strictly acylindrical colorable HHG is
noetherian.
\end{restatableThm}

The colorability assumption in \Thmref{main-thm} (see \defref{colorable-HHG})
is not a restricting one, and is supported by numerous examples, including
$\left(i\right)$ mapping class group \cite{BBF}; $\left(ii\right)$
many cubical groups (including all right-angled Artin and Coxeter
groups); \cite{BHS1,Ha20} $\left(iii\right)$ Teichmüller space (with
either metrics) \cite{Ra07,Dur,EMR}, among others. An example of
non-colorable HHG has recently been provided in \cite{Ha23}.

In contrast, the assumption of strict acylindricity (see \defref{strictly-acylindrical-product}
and \defref{stricly-acylidrical-HHG}) is a limiting one, and does
not hold in the case of mapping class group. This assumption is a
strengthening of weak acylindricity (see \cite[Definition 2.2]{Se23}),
which does holds for the mapping class group. We hope that the ideas
in this paper could be generalized to weakly acylindrical groups.
\begin{rem}
We suspect that as a possible application of \Thmref{main-thm}, which
we will not pursue in this paper, one can construct a canonical higher-rank
version of the Makanin-Razborov diagram, generalizing the canonical
Makanin-Razborov diagram for hyperbolic groups \cite{Se1,We-Re,Se7}
as well as the non-canonical higher-rank diagram of colorable HHGs
(\cite{Se23,JA-SE}).
\end{rem}

The paper is organized as follows: In \Secref{Acylindrical-groups}
we provide an overview of limit groups over acylindrically hyperbolic
groups, primarily based on \cite{GR-HU}, while also incorporating
certain proofs and definitions that are not explicitly stated but
are essential for this work. In \Secref{Case-of-Product} we prove
the main theorem in a specific case involving product of hyperbolic
spaces, aiming to illustrate the proof's concepts. Finally, in \Secref{general-digram},
we utilize the concept of projection complexes and quasi-tree of metric
spaces from \cite{BBF} and \cite{BBFS}. These concepts enable us
to apply the notion of limit groups to colorable HHGs and subsequently
prove the main theorem.
\begin{rem}
For clarity, this paper presupposes familiarity with HHGs, limit groups,
and JSJ. We reiterate essential definitions and properties while directing
readers to external sources as necessary.
\end{rem}

\section{\protect\label{sec:Acylindrical-groups}Acylindrical groups - and
Limit Groups over them}

\global\long\def\im{\text{Im}}%
\global\long\def\conj{\text{Conj}}%
\global\long\def\res{\text{Res}}%
\global\long\def\acts{\curvearrowright}%
\global\long\def\jsj{\text{JSJ}}%
\global\long\def\surject{\twoheadrightarrow}%
\global\long\def\comp{\text{Comp}}%
\global\long\def\sker{\underrightarrow{\ker}^{\omega}}%
\global\long\def\amal#1#2#3{#1\underset{#3}{*}#2}%
 
\global\long\def\hom{\text{Hom}}%
\global\long\def\isom{\text{Isom}}%
\global\long\def\bndv{\text{qTriv}}%
\global\long\def\bnd{\bndv_{c}}%
\global\long\def\wsurelu{\ \omega\text{-almost surely}}%
\global\long\def\onto{\twoheadrightarrow}%
\global\long\def\restriction#1#2{\left.\kern-\nulldelimiterspace#1\vphantom{\big|}\right|_{#2}}%
\global\long\def\mod{\text{Mod}}%
\global\long\def\aut{\text{Aut}}%
\global\long\def\stab{\text{Stab}}%
\global\long\def\by{\mathbf{Y}}%
 
\global\long\def\projcopx{\mathcal{P}_{K}}%
\global\long\def\diam{\text{diam}}%
\global\long\def\qtms{\mathcal{C}_{K}}%
\global\long\def\c{\mathcal{C}}%
\global\long\def\ufin{\mathcal{A}_{\text{ufin}}}%
\global\long\def\vab{\mathcal{A}_{\infty}}%
\global\long\def\nvc{\mathcal{A}_{\mathrm{nvc}}}%
\global\long\def\ahtree{\left(\mathcal{A},\mathcal{H}\right)\text{-tree}}%
\global\long\def\ah{\left(\mathcal{A},\mathcal{H}\right)}%
\global\long\def\ahjsj{\ah\text{-JSJ}}%
\global\long\def\qh{QH}%
\global\long\def\inc{\text{Inc}_{v}^{\mathcal{H}}}%
 
\global\long\def\jsj{\mathbb{A}_{JSJ}}%
\global\long\def\mjsj{\mathbb{A}_{M}}%
\global\long\def\linnell#1{\mathbb{D}_{#1}}%
\global\long\def\lin{\linnell L}%
\global\long\def\sbt{\mathcal{SBT}}%
\global\long\def\lsbt{\sbt}%
\global\long\def\lse{\mathcal{LSE}}%

In this section we give a quick review of the topic of limit groups
over acylindrically hyperbolic groups. We will mostly skim the results
from \cite{GR-HU}, with slight modification, mentioned when needed.

For the rest of the paper we fix some non-principal ultrafilter $\omega$
on $\mathbb{N}$.

\begin{defn}
Let $\Gamma$ be a group, $G$ a finitely generated group, and $\left(\varphi_{n}\right)_{n}$
a sequence in $\hom\left(G,\Gamma\right)$. The \textbf{$\omega$-stable
kernel} of $\left(\varphi_{n}\right)$ is
\[
\sker\varphi_{n}=\left\{ g\in G\mid\varphi_{n}\left(g\right)=1\ \wsurelu\right\} 
\]

The associated \textbf{$\Gamma$-limit group} (or just limit group,
if $\Gamma$ is clear) is the quotient $L=\nicefrac{G}{\sker\varphi_{n}}$.
We call the map $\varphi_{\infty}:G\onto L$ the \textbf{limit quotient}
map and $\left(\varphi_{n}\right)$ the \textbf{defining sequence}
of $L$.
\end{defn}

\begin{rem}
Note the same limit group might be obtained by many different sequences,
i.e. a ``defining sequence'' is not unique. For the majority of
the paper, the definition of limit group includes the defining sequence.
In the definition of JSJ decomposition of limit groups (see \subsecref{JSJ-decomposition})
we take into account all possible defining sequences.
\end{rem}

In our context, limit groups relate to equational noetherianity by
the following easy lemma.
\begin{lem}
\label{lem:en_is_eq_factoring}\cite[Lemma 3.5]{GR-HU}For a group
$\Gamma$, the following are equivalent:
\end{lem}

\begin{enumerate}
\item $\Gamma$ is equationally noetherian.
\item For any finitely generated group $G$ and any sequence $\left(\varphi_{n}\right)\in\hom\left(G,\Gamma\right)$
we have $\varphi_{n}$ $\wsurelu$ factors through the limit quotient
$\varphi_{\infty}:G\to L$.
\item For any finitely generated group $G$ and any sequence $\left(\varphi_{n}\right)\in\hom\left(G,\Gamma\right)$
there exists $n$ such that $\varphi_{n}$ factors through the limit
quotient $\varphi_{\infty}:G\to L$.
\end{enumerate}
Using this lemma is the preferred strategy to prove that certain groups
are equationally noetherian \cite{GR-HU,Se7,Se10,We-Re}, and this
is the strategy used in this paper. In addition to \ref{lem:en_is_eq_factoring},
we will need the following lemma, which we will use in proving \ref{prop:step-induction-product}.
The proof is not difficult, and can be found in \cite[Lemnma 3.19]{GR-HU}.
\begin{lem}[{\cite[Lemnma 3.19]{GR-HU}}]
\label{lem:relative_presented_lemma}Suppose that $\left(\varphi_{n}\right)\in\hom\left(G,\Gamma\right)$
is such that the limit group $L=\nicefrac{G}{\sker\left(\varphi_{n}\right)}$
is finitely presented relative to subgroups $\left\{ P_{1},\ldots,P_{n}\right\} $.
Suppose further that for each $P_{j}$, there is a subgroup $\tilde{P}_{j}\leq G$
which $\varphi_{\infty}$ maps onto $P_{j}$ such that $\restriction{\varphi_{n}}{\tilde{P}_{j}}$
factors through $\restriction{\varphi_{\infty}}{\tilde{P}_{j}}$ $\wsurelu$.
Then $\varphi_{n}$ factors through $\varphi_{\infty}$ $\wsurelu$.
\end{lem}

\begin{defn}
Let $\Gamma$ be a group acting by isometries on a $\delta$-hyperbolic
metric space $X$.

The action is called \textbf{acylindrical} if for every $\varepsilon>0$
there exist $R,N>0$ such that for any $x,y\in X$ such that $d\left(x,y\right)>R$,
the set
\[
\left\{ g\in\Gamma\mid d\left(x,g\cdot x\right),d\left(y,g\cdot y\right)<\varepsilon\right\} 
\]

is of cardinality at most $N$. This property can be thought of as
``coarse pairwise stabilizers are uniformly finite''.
\end{defn}

\begin{defn}
Let $\Gamma$ be a group acting by isometries on a $\delta$-hyperbolic
metric space $X$. The action of $\Gamma$ on $X$ is called \textbf{elliptic}
if the set $\Gamma\cdot x$ is bounded for any, or equivalently some,
$x\in X$.

The action is called non-elementary if $\Gamma$ is not virtually
cyclic and not elliptic. If $\Gamma$ admits a non-elementary acylindrical
action on some $\delta$-hyperbolic $X$ we say that $\Gamma$ is
\textbf{acylindrically hyperbolic }and call $X$ an auxiliary space.
\end{defn}

\begin{defn}[{\cite[Definition 4.2]{GR-HU}}]
{\footnotesize\label{def:divergent_seq}}Let $\Gamma$ be a group
acting on a metric space $X$ and let $G$ be a group generated by
a finite set $S$. Let $\varphi:G\to\Gamma$ be a homomorphism. Define
the \textbf{$X$-scaling factor} of $\varphi$ by
\[
\left\Vert \varphi\right\Vert _{X}=\inf_{x\in X}\max_{s\in S}d\left(x,\varphi\left(s\right)\cdot x\right).
\]
If the context is clear we denote $\left\Vert \varphi\right\Vert =\left\Vert \varphi\right\Vert _{X}$
and call this the \textbf{scaling factor} of $\varphi$. If $\left(\varphi_{n}\right)$
is a sequence of homomorphisms from $G$ to $\Gamma$ the we say that
$\left(\varphi_{n}\right)$ is \textbf{$X$-divergent }(or just \textbf{divergent})
if
\[
\lim_{\omega}\left\Vert \varphi_{n}\right\Vert =\infty
\]
and \textbf{non-divergent }if $\lim_{\omega}\left\Vert \varphi_{n}\right\Vert <\infty$.
We say that the limit group corresponding to $\left(\varphi_{n}\right)$
divergent or non-divergent (with the defining sequence $\left(\varphi_{n}\right)$)
correspondingly.

{\tiny}{\tiny\par}
\end{defn}

\subsection{\protect\label{subsec:JSJ-decomposition}JSJ Decomposition of Limit
Groups over Acylindrical Action and Modular Automorphisms}

Fix some \textbf{acylindrically hyperbolic} group $\Gamma$ with auxiliary
$\delta$-hyperbolic space $X$. In this section we will describe
the JSJ decomposition of divergent limit groups, without providing
many proofs.

In order to use the language of Guirardel-Levitt \cite{Gu-Le}, we
need to recall a few notions relating trees with group actions. All
of these notions appear in \cite{GR-HU}.
\begin{defn}
Given a group $G$, let $\mathcal{A}$ and $\mathcal{H}$ be two families
of subgroups such that $\mathcal{A}$ is closed under conjugation
and taking subgroups.
\begin{itemize}
\item An \textbf{$\bm{\ahtree}$} is a simplicial tree $T$ together with
an action of $G$ on T such that each edge stabilizer belongs to $\mathcal{A}$
and each $H\in\mathcal{H}$ fixes a vertex in $T$. An \textbf{$\left(\mathcal{A},\mathcal{H}\right)$-splitting
}is a graph of groups decomposition whose corresponding Bass-Serre
tree is an $\left(\mathcal{A},\mathcal{H}\right)$-tree. This is also
referred to as a splitting over $\mathcal{A}$ relative to $\mathcal{H}$
(or a tree over $\mathcal{A}$ relative to $\mathcal{H}$).
\item An $(\mathcal{A},\mathcal{H})$-tree $T$ is called \textbf{universally
elliptic} if the edge stabilizers of $T$ are elliptic in every other
$(\mathcal{A},\mathcal{H})$-tree. Also, given trees $T$ and $T^{\prime}$,
we say $T$ \textbf{dominates} $T'$ if every vertex stabilizer of
$T$ is elliptic with respect to $T'$.
\item Given an $\ahtree$ $T$ and an edge $e$ in $T$, one can collapse
the edge $e$ in an equivariantly to get a new tree. A \textbf{collapse
map} is a finite composition of such collapses. If we have a collapse
map $T\to T'$ we say that $T$ is a \textbf{refinement} of $T'$.
Two trees are \textbf{compatible} if they have a common refinement.
\item An \textbf{$\bm{\ahtree}$} $T$ of $G$ is called an \textbf{$\bm{(\mathcal{A},\mathcal{H})}$--JSJ}
tree if $T$ is an $\ahtree$ which is universally elliptic and which
dominates every other universally elliptic $\ahtree$. In this case
the associated $\left(\mathcal{A},\mathcal{H}\right)$-splitting is
called an $\left(\mathcal{A},\mathcal{H}\right)$-JSJ decomposition.
\item Given a vertex $v$ in an $\ahtree$ $T$ we denote by $\bm{\inc}$
the set of subgroups of the vertex group $G_{v}$ which are conjugate
into either the stabilizer of some edge group of an edge adjacent
to $v$ or into some group contained in $H\in\mathcal{H}$. 
\item A vertex group $Q=G_{v}$ of an $(\mathcal{A},\mathcal{H})$--splitting
of $G$ is called a $\bm{\qh}$\textbf{-subgroup} if there is a normal
subgroup $N\trianglelefteq Q$ (called the fiber of $Q$) and a hyperbolic
$2$-orbifold $\Sigma$ such that the following hold:
\begin{itemize}
\item $\nicefrac{Q}{N}$ is isomorphic to the fundamental group $\Sigma$,
$\pi_{1}\left(\Sigma\right)$, and
\item if $e$ is an edge adjacent to $v$, the image of $G_{e}$ in $\nicefrac{Q}{N}$
is either finite or conjugate to the fundamental group of a boundary
component of the orbifold.
\end{itemize}
\item A $\qh_{(\mathcal{A},\mathcal{H})}$-subgroup is a $\qh$-subgroup
$Q$ such that: 
\begin{itemize}
\item For every essential simple closed curve $\gamma$ on the underlying
orbifold $\Sigma$, the corresponding subgroup $Q_{\gamma}$ of $\pi_{1}\left(\Sigma\right)=\nicefrac{Q}{N}$
is the image of a subgroup in $\mathcal{A}$ under the quotient map
$Q\twoheadrightarrow\nicefrac{Q}{N}$. Further,
\item the intersection of $Q$ and any element of $H\in\mathcal{H}$ has
image in $\nicefrac{Q}{N}$ that is either finite or conjugate to
the fundamental group of a boundary component of the orbifold.
\end{itemize}
\begin{rem*}
The isomorphism type of $Q$ does not necessarily determine the surface
nor the fiber. When we refer to $\qh_{\left(\mathcal{A},\mathcal{H}\right)}$
subgroup, we always consider $N$ and $\Sigma$ as part of the structure.
\end{rem*}
\item A group is called $C$-virtually cyclic if it maps onto either $\mathbb{Z}$
or $D_{\infty}$ with kernel of order at most $C$.
\end{itemize}
\end{defn}

We fix some divergent limit group $L$ over $\Gamma$ with defining
sequence $\left(\varphi_{n}\right)$.

we take a slightly different route from \cite{GR-HU}. Our approach
is not essentially different, but allows us a bit more control over
descending chains of limit groups (see \subsecref{Descending-Chains-Condition}).
We will note when diverging from the definitions in \cite{GR-HU}.
\begin{defn}
A finitely generated subgroup $H$ of $L$ is the called \textbf{stably
elliptic} if for some finitely generated $\varphi_{\infty}$-lift
$\tilde{H}\leq G$, the subgroups $\varphi_{n}\left(\tilde{H}_{n}\right)$
of $\Gamma$ are $\wsurelu$ elliptic. A subgroup $H$ of $L$ is
called \textbf{locally stably elliptic} if every finitely generated
subgroup $H_{0}\leq H$ is stably elliptic, denote the set of locally
stably elliptic subgroup of $L$ with respect to $\left(\varphi_{n}\right)$
by $\lse\left(L,\left(\varphi_{n}\right)\right)$.
\end{defn}

\begin{rem}
Note that those definition do not depend on the choice of lifts, but
does depend on the defining sequence of $L$.
\end{rem}

If $H\in\lse\left(L,\left(\varphi_{n}\right)\right)$ then in the
action of $L$ on the associated limit $\mathbb{R}$-tree, $H$ fixes
a point (\cite[Lemma 4.7, 5]{GR-HU}) thus using this action one can
not hope to learn about the algebraic structure of such groups (indeed,
those subgroups could be arbitrary, e.g. for any group $K$, the group
$\mathbb{Z}*K$ is acylindrically hyperbolic with $K$ elliptic).
Despite of that, if $H$ is \textbf{not }locally stably elliptic with
respect to some other defining sequence (i.e. with respect to $\left(\psi_{n}\right)\in\hom\left(G',\Gamma\right)$
such that $\nicefrac{G'}{\sker\psi_{n}}=L$) , meaning $H$ ``see''
the hyperbolicity of $\Gamma$ is some way, then we still can use
the other defining sequence to learn about $H$. For more details
see \cite[Section 5]{GR-HU}. This motivate the definition of absolutely
elliptic subgroup of $L$, which are never ``sees'' the hyperbolicity
of $\Gamma$.
\begin{defn}[{\cite[Def 5.4]{GR-HU}}]
 We say that $H\leq L$ is \textbf{absolutely elliptic }if for any
finitely generated group $G_{0}$ and any sequence $\left(\psi_{n}\right)\in\hom\left(G_{0},\Gamma\right)$
that define $L$ (i.e. $L=\nicefrac{G_{0}}{\sker\psi_{n}}$) we have
$H\in\lse\left(L,\left(\psi_{n}\right)\right)$.
\end{defn}

In \cite{GR-HU} the authors build the JSJ decomposition of a divergent
limit group using the notion of absolutely elliptic subgroups. In
out context, in order to prove a descending chain condition, we need
a bit more control over subgroups of $L$ that are by themselves ``non-divergent''.
\begin{defn}
\label{def:boundedness}An element $g\in L$ is called \textbf{stably}
\textbf{bounded-trace} if for some (thus all) lifts $\tilde{g}\in G$
of $g$ we have $\lim_{\omega}\inf_{x\in X}d\left(x,\varphi_{n}\left(\tilde{g}\right)x\right)<\infty$.
Denote the set of \textbf{stably bounded-trace elements} by $B_{L}$.

A finitely generated subgroup $H\leq L$ is called \textbf{stably
bounded-trace} if for some (thus all) finitely generated lift $\tilde{H}\le G$
we have $\lim_{\omega}\left\Vert \restriction{\varphi_{n}}{\tilde{H}}\right\Vert =\lim_{\omega}\inf_{x\in X}\max_{s\in S}d\left(x,\varphi_{n}\left(s\right)\right)<\infty$
where is $S$ is some finite generating set of $\tilde{H}$. Denote
the set of locally stably bounded trace, with respect to $\left(\varphi_{n}\right)$
by $\lsbt\left(L,\left(\varphi_{n}\right)\right)$.

This condition on $H_{0}$ can be read as ``if we restrict $\left(\varphi_{n}\right)$
to $\tilde{H}_{0}$ we get a non-divergent sequence''. Every locally
stably elliptic subgroup is locally stably bounded trace subgroup
(\cite[Lemma 2.2]{GR-HU} gives a uniform bound on the minimal displacement
of stable elliptic subgroups of a limit group).

\end{defn}

\begin{rem}
Note that $B_{L}$ is generally not a subgroup, but it is closed under
conjugation and taking powers. If $L$ is a non-divergent limit group
then $B_{L}=L$.
\end{rem}

In the following definition we first diverge a bit from \cite{GR-HU}.
Out approach is not essentially different, and the proof of existence
and the structure of the JSJ (\thmref{JSJ-exists}) remains valid
since we only enlarged the set of subgroups we demand to be elliptic
in the JSJ.
\begin{defn}[{\cite[Definitions 4.3 \& 5.12]{GR-HU}}]
Fix $C>0$.

Denote by $\ufin$ the set all finite subgroup of $L$ of order at
most $2C$, $\vab$ all virtually abelian subgroups of $L$ which
are not absolutely-elliptic. Let $\mathcal{A}=\vab\cup\ufin$. We
also denote $\nvc$ the set of groups in $\mathcal{A}$ that are not
virtually cyclic.

Let $\mathcal{H}$ be the set of finitely generated stable bounded-trace
subgroups of $L$.
\end{defn}

\begin{rem}
The reason behind using the notion of absolutely-elliptic (and absolutely
bounded trace) subgroups is that the JSJ of a limit group should not
depend on its defining sequence. A fact that is used in the in whats
called shortening argument. In our case, in order the prove descending
chain condition (see \ref{subsec:Descending-Chains-Condition}), we
do need all locally bounded trace bounded to be elliptic in the JSJ
(and not ``absolutely bounded trace''). This difference does not
change the results of the shortening argument in \cite{GR-HU}, in
\ref{subsec:Modular-automorphism} we note on this difference.
\end{rem}

\begin{thm}[{\cite[Theorem 5.15 \& Lemma 5.17]{GR-HU}}]
\label{thm:JSJ-exists}Let $L$ be a divergent limit group that does
not split over $\ufin$. Then there exists a $\ahjsj$ decomposition
$T$ of $L$ such that every vertex $v$ is either:
\begin{enumerate}
\item rigid, i.e. $L_{v}$ is elliptic in any $\ah$-tree
\item flexible (i.e. not rigid), and then either:
\begin{enumerate}
\item $L_{v}$ is virtually abelian
\item \label{enu:flexible-vertex}$L_{v}$ admit a splitting over $C$-virtually
cyclic absolutely elliptic subgroups and relative to $\inc$ in which
vertex groups that are each either $\left(i\right)$ rigid in any
$\ah$-tree or $\left(ii\right)$ $\qh_{\ah}$ with finite fiber.
\end{enumerate}
\end{enumerate}
\end{thm}

\begin{rem}
By definition of JSJ decompositions, every $\ah$ one-edge splitting
of $L$ is either a collapse of the JSJ decomposition, or comes from
splitting one of the flexible vertex groups of the JSJ decomposition,
and then collapsing all other edges.

Let $v$ be a flexible vertex of type \ref{enu:flexible-vertex}.
This condition means that any refinement of $L_{v}$ relative to $\inc$
is done by splitting some $\qh_{\ah}$ vertex groups along some simple
closed curves and then collapsing some edges.

\end{rem}

\begin{prop}[{\cite[Theorem 5.20]{GR-HU}}]
\label{prop:tree-of-cylinders}Given a divergent limit group $L$
that does not split over $\ufin$, and $T$ an $\ahjsj$ of $L$,
there exists an $\ah$-tree $T_{c}$ called the tree of cylinders
such that
\begin{enumerate}
\item $T_{c}$ is $\left(2,C\right)$-acylindrical (i.e. every point wise
stabilizer of a path of length $\geq3$ has order $\leq C$).
\item $T$ dominates $T_{c}$
\item Any group which is elliptic in $T_{c}$ but not in $T$ is virtually
abelian and not virtually cyclic.
\end{enumerate}
Furthermore, $T_{c}$ is a $\left(\mathcal{A},\mathcal{H}\cup\nvc\right)\text{-JSJ}$
tree and is compatible with every $\ah$-tree.

\end{prop}

\begin{defn}[{\cite[Def 5.21]{GR-HU}}]
 Let $L$ be a non-divergent limit group that does not split over
$\ufin$. The \textbf{JSJ tree} and \textbf{JSJ decomposition} of
$L$ (with out specification of families of groups) are $T_{c}$ from
\Propref{tree-of-cylinders} and its corresponding graph of groups,
denoted $\mathbb{A}_{JSJ}$. The \textbf{modular tree}, denoted\textbf{
}$T_{M}$, is obtained by refining $T_{c}$ at every flexible vertex
group satisfying \enuref{flexible-vertex} in \Thmref{JSJ-exists},
and the corresponding decomposition is called the \textbf{modular
splitting }of $L$ and denoted $\mjsj$.
\end{defn}

\begin{rem}
Vertex groups in $\mathbb{A}_{M}$ can be characterized as $\qh_{\ah}$,
virtually abelian and rigid in the usual sense.
\end{rem}

Now assume $L$ that does split over $\ufin$ (or might split over
$\ufin$). In this case we look at $\lin$, the Linnell Decomposition
of $L$, which is a $\left(\ufin,\mathcal{H}\right)$-JSJ decomposition
in the sense of Guirardel-Levitt (see \cite{Li83}). In this decomposition,
every edge group is finite and every vertex group is either finite
or does not splits over $\ufin$ relative to $\mathcal{H}$.

Let $L_{v}$ be some vertex group in $\mathbb{D}_{L}$. The vertex
group $L_{v}$ is finitely generated, thus its also a $\Gamma$-limit
group at its own right (restrict $\left(\varphi_{n}\right)$ to a
finitely generated lift of $L_{v}$). We will associate a splitting
to $L_{v}$:
\begin{itemize}
\item If $L_{v}$ is not absolutely bounded trace subgroup then has a $\left(\mathcal{A}_{\infty},\mathcal{H}\cup\nvc\right)\text{-JSJ}$
splitting (as one can realize $L_{v}$ as a non-divergent limit group,
and thus it has a JSJ), denoted $\mathbb{A}_{v}$.
\item If $L_{v}$ is absolutely bounded trace subgroup or finite then we
take the trivial decomposition as $\mathbb{A}_{v}$.
\end{itemize}
Thanks to the fact that finite subgroups are elliptic in any tree,
we can glue all the splittings $\mathbb{A}_{v}$ along the edges of
$\mathbb{D}_{L}$ to get a decomposition of $L$. One can think of
this splitting as a blowup of $\mathbb{D}_{L}$ at every vertex which
has non absolutely bounded vertex stabilizer.
\begin{defn}
\label{def:ufin_JSJ}Abusing some notation, we also call this resulting
splitting the JSJ of $L$ and denote it by $T_{JSJ}$ and the resulting
decomposition by $\mathbb{A}_{JSJ}$ (even though it in not a JSJ
in the sense of Guirardel-Levitt). We consider the vertices in $T_{JSJ}$
that correspond to vertices in $\mathbb{D}_{L}$ with finite stabilizer
or with stabilizers (i.e. the one in the that absolutely bounded-trace
as rigid vertex subgroup of $T_{JSJ}$\textbackslash$\jsj$. The
rest of vertices of $T_{JSJ}$ are colored as rigid/virtually abelian/QH
as their corresponding coloring in the JSJ decomposition of the vertex
stabilizer they correspond to.
\end{defn}

\subsection{Modular Automorphism\protect\label{subsec:Modular-automorphism}}

Given a divergent limit group $L$ , keeping the assumption that $L$
does not splits over $\ufin$, one can define the modular automorphisms
group $\mod\left(L\right)$ to be essentially the subgroup of $\aut\left(L\right)$
generated by inner automorphisms, Dehn twists over edges in the modular
splitting, and extension of automorphisms of flexible vertex group
$v$ respecting the edge groups adjacent to $v$. See \cite[Section 5.4]{GR-HU}
for the full definition.

We continue as in \cite[Section 6]{GR-HU} and \cite[section 7]{We-Re},
and associate to $L$ a series of finitely presented covers $\left\{ W_{i}\right\} $
converging to $L$ (i.e. $L$ is the direct limit of $\left\{ W_{i}\right\} $),
such that each of the $W_{i}$ has a decomposition similar to the
JSJ of $L$ (the decomposition has the same underlying graph as $\mathbb{A}_{JSJ}$,
and it satisfies more conditions, see \cite[Lemma 6.3]{GR-HU} and
\cite[Lemma 7.1]{We-Re}). Since each $W_{i}$ is finitely presented
and has $L$ as quotient then $\wsurelu$ $\left(\varphi_{n}\right)$
factor through $W_{i}$. We can move to subsequences of $\left(\varphi_{n}\right)$
and re-indexing we can assume that if $i>n$ then $\varphi_{n}$ factor
through $W_{i}$.

Denote by $\hat{\varphi}_{n}:W_{n}\to\Gamma$ the map associated to
$\varphi_{n}:G\to\Gamma$, i.e the map we get from the fact that $\varphi_{n}$
factor through $W_{n}$. Among all maps $W_{n}\to\Gamma$ of the form
$\hat{\varphi}_{n}\circ\alpha$ where $\alpha\in\mod\left(W_{n}\right)$
(defined the same way as $\mod\left(L\right)$) choose the one minimizing
the scaling factor (as maps to $\Gamma$), i.e. choose $\alpha\in\mod\left(W_{n}\right)$
such that for any $\beta\in\mod\left(W_{i}\right)$ we have $\left\Vert \hat{\varphi}_{n}\circ\alpha\right\Vert \leq\left\Vert \hat{\varphi}_{n}\circ\beta\right\Vert $
(this $\alpha$ is not unique). Doing this for every $n$ we get a
new sequence from $\psi_{n}:=\hat{\varphi}_{n}\circ\alpha_{n}\circ\xi_{n}:G\to\Gamma$
(where $\xi_{n}:G\to W_{n}$ is the quotient map which comes from
the definition of $W_{n}$). This sequence defines a new limit group
$L_{2}$ which is a quotient of $L$. We call $L_{2}$ a \textbf{shortening
quotient} of $L_{2}$.

\paragraph*{For a limit group $L$ which does split over $\protect\ufin$}

one can build similar sequence $\left\{ W_{n}\right\} $ as in the
case as in the case where $L$ does not split over $\ufin$ (take
the series corresponding of finitely presented coverts of each $L_{v}$
and glue them as in the construction of \ref{def:ufin_JSJ}) Now define
the shortening quotient as before. In \cite{GR-HU} and \cite{We-Re}
they formulate this construction in a different way, but its essentially
identical.

The following fact, which we won't proof, is a standard fact about
shortening quotient followed easily from its definition (for details
see \cite[Section 6]{GR-HU} and especially the proof of \cite[Lemma 6.5]{GR-HU}).
\begin{lem}
\label{lem:shortening_on_rigid_is_conj}Keeping the notations of the
discussion above, let $\psi_{\infty}:G\to L_{2}$ be the limit quotient
of $\left(\psi_{n}\right)$ and denote by $\pi:L\to L_{2}$ the shortening
quotient map. Let $V$ be a rigid vertex group of the JSJ decomposition
of $L$ and $V_{0}\leq V$ a finitely generated subgroup of $V$,
and $\tilde{V}_{0}\leq G$ a $\varphi_{\infty}$-lift of $V_{0}$.
Then from some point onward we have that $\restriction{\varphi_{n}}{V_{0}}$
and $\restriction{\psi_{n}}{\tilde{V_{0}}}$ differ by post conjugation
in $\Gamma$.
\end{lem}

The proof of the following lemma from \cite{GR-HU} is essentially
the shortening argument of Sela, modified.
\begin{lem}[{\cite[Lemma 6.4 and Lemma 6.6]{GR-HU}}]
\label{lem:shortening-argument}Keeping the notations of the discussion
above, let $\psi_{\infty}:G\to L_{2}$ be the limit quotient of $\left(\psi_{n}\right)$
and denote by $\pi:L\to L_{2}$ the shortening quotient map. Then:
\begin{enumerate}
\item $\pi$ is injective restricted to rigid vertex groups of $L$.
\item If $\pi$ is an isomorphism then $L_{2}$ is a non-divergent limit
group by itself (i.e. $\lim_{\omega}\left\Vert \psi_{i}\right\Vert <\infty$).
\end{enumerate}
\end{lem}

\begin{proof}
Clause $\left(1\right)$ follows in the same way as in \cite{GR-HU}.
Clause $\left(2\right)$ also follows in a similar way, but one need
to note few facts beforehand. Assume that $\pi$ is an isomorphism,
and that $\left(\psi_{n}\right)$ is non-divergent. Let $T$ be the
$\mathbb{R}$-tree corresponding to $\left(\psi_{n}\right)$, and
let $\mathbb{A}$ the the splitting obtained by the Rips machine.
The proof of for the shortening quotient \cite[Theorem 5.29]{GR-HU}
uses the fact that $\mathbb{A}$ is a $\left(\mathcal{A},\mathcal{H}_{0}\right)$-tree
where $\mathcal{H}_{0}$ is the set of all finitely generated absolutely
elliptic subgroups, and thus $\mod\left(\mathbb{A}\right)\leq\mod\left(L\right)$
(\cite[5.23]{GR-HU}). This is true in out case also, since we did
not change the set of allowed edge stabilizers and all bounded trace
subgroups are also elliptic in $T$ (the proof of \cite[Lemma 4.7 (7)]{GR-HU}
works also for finitely generated locally bounded trace subgroups).
Now one can apply the same shortening argument.
\end{proof}
In light of the definition of shortening quotient, given a limit group
$L_{1}$, one can short (i.e. build shortening quotient) $L_{1}$
to get the shortening quotient $L_{1}\onto L_{2}$, and proceed repetitively
to get a sequence $L_{1}\onto L_{2}\onto L_{3}\onto\cdots$. A priori,
this procedure might continue indefinitely. Notice that every shortening
quotient commute with the limit quotients of the corresponding limit
groups, motivating the following definition
\begin{defn}
A \textbf{resolution} is a sequence of limit groups with corresponding
surjective maps that commute with the limits quotients, i.e. $L_{1}\overset{\tau_{1}}{\onto}L_{2}\overset{\tau_{2}}{\onto}\cdots\overset{\tau_{n-1}}{\onto}L_{n}\onto\cdots$
where $L_{i}$ is a limit group with with corresponding limit quotient
$\varphi_{\infty}^{i}:G\onto L_{i}$ such that $\varphi_{\infty}^{i}\circ\tau_{i}=\varphi_{\infty}^{i+1}$.
The sequence might be finite or infinite.
\end{defn}

\subsection{\protect\label{subsec:Totally-rigid-subgroup}Totally Rigid Subgroup}

In this section we explicitly discuss an idea appearing in end of
the proof of \cite[Proposition 13]{JA-SE} (for the enthusiastic reader,
the mentioned idea appears at the end of the proof of Proposition
17).

Let $L_{1}$ be a limit group that does not admit an infinite resolution
\[
L_{1}\onto L_{2}\onto\cdots
\]

In \cite[Lemma 6.7]{GR-HU}, the authors proves that a limit group
has an infinitely generated subgroup which is virtually abelian and
non-(absolutely elliptic) then the limit group admits an infinite
sequence of proper limit quotients. This works in our setting too,
since the addition of absolutely bounded-trace groups as elliptic
in the JSJ does not change any of their arguments.

This yields the following lemma:
\begin{lem}
\label{lem:non-(absolutly-elliptic)-vAb-are-f.g.}All non-(absolutely
elliptic) virtually abelian subgroups of $L_{1}$ and of any limit
group which is quotient of $L_{1}$ are finitely generated. In particular,
all the edge groups in the JSJ decompositions of said limit groups
are finitely generated.
\end{lem}

Let $L_{1}\onto L_{2}\onto\cdots\onto L_{s}$ be some resolution obtained
by successive shortening quotients of $L_{1}$.

\paragraph*{Our Goal:}

To find finitely many finitely generated subgroups $P_{1},\ldots,P_{m}$
of $L_{1}$ such that a) $L_{1}$ is finitely presented relative to
$P_{1},\ldots,P_{m}$ and b) the image of each $P_{1},\ldots,P_{s}$
in every group in the resolution (including $L_{1}$) is contained
in some rigid vertex group (which might be different for each $P_{i}$)
of the the JSJ of the corresponding group (hence ``totally rigid'').

If we find such $P_{1},\ldots,P_{m}$, then since shortening quotients
are injective when restricted to rigid vertex groups, all the $P_{i}$
will be contained in all the limit groups $L_{i}$ in the resolution.

For a start, let us look at the quotient $L_{1}\onto L_{2}$. Denote
by $T_{M}^{1}$ and $T_{M}^{2}$ and $\mathbb{A}_{M}^{1}$ and $\mathbb{A}_{M}^{2}$
the modular JSJ trees and decomposition of $L_{1}$ and $L_{2}$ correspondingly.

Let $P$ be some finitely generated subgroup of $L_{1}$ contained
in some rigid vertex group $T_{M}^{1}$. The group $P$ is also a
subgroup of $L_{2}$ (since shortening quotient are injective in restriction
to rigid vertex groups), thus $P$ acts on $T_{M}^{2}$. Denote the
$P$-minimal subtree of $T_{M}^{2}$ by $T_{2}$. Denote the graph
of groups associated to the action $P$ of $P$ on $T_{2}$ by $\Delta_{2}$.

The vertices of $T_{M}^{2}$ are colored as rigid, virtually abelian
and $\qh$, and we can restrict this coloring to $T_{2}$, and $\Delta_{2}$
(notice that the coloring of $T_{M}^{2}$ is equivariant, thus so
is the coloring of $T_{2}$).

Note that $\Delta_{2}$ is a finite graph of groups (since $P$ is
finitely generated and $T_{2}$ is minimal). Denote by $Q_{1},\ldots,Q_{k}$
the rigid subgroups in $\Delta_{2}$.
\begin{prop}
\label{prop:totally-rigid-step-1}$P$ is finitely presented relative
to $Q_{1},\ldots,Q_{k}$.
\end{prop}

The proof will be divided into cases, showing that non-rigid vertex
groups of $\Delta_{2}$ must be finitely presented, and edge groups
of $\Delta_{2}$ must be finitely generated.
\begin{lem}
Let $E$ be an edge group in $T_{2}$, then $E$ is finitely generated
\end{lem}

\begin{proof}
Denote by $e$ the edge stabilized by $E$.

We have three options for $e$:
\begin{enumerate}
\item The edge $e$ comes from an edge in the Linnell Decomposition of $L_{2}$.
\item The edge $e$ comes from the splitting of flexible subgroup done while
building $\mathbb{A}_{M}^{2}$ from $\mathbb{A}_{JSJ}^{2}$, the JSJ
decomposition of $L_{2}$ (i.e. $e$ comes from a blowup of some flexible
splitting from \Thmref{JSJ-exists} satisfying \ref{enu:flexible-vertex}).
\item The edge $e$ comes from $\mathbb{A}_{JSJ}^{2}$.
\end{enumerate}
In cases 1 and 2, $E$ in contained either in a uniformly finite subgroup
of $L_{2}$ or in a $C$-virtually cyclic group correspondingly. In
both cases $E$ is finitely generated. In case $3$ the group $E$
is contained in a finitely generated virtually abelian group, thus
also finitely generated.
\end{proof}
\begin{lem}
All vertex groups of $T_{2}$ are finitely generated
\end{lem}

\begin{proof}
If we quotient $P$ by all edge groups in $\Delta_{2}$ we get the
free product of the vertex groups quotiented by their edge groups.

Since $P$ is finitely generated, the factors in the free product
are also finitely generated. Each vertex group is finitely generated
relative to its edge groups. Since the edge groups are finitely generated
as proved before, it follows that the vertex groups are themself finitely
generated.
\end{proof}
\begin{lem}
Let $V$ be a $\qh$ vertex group in $T_{2}$, then $V$ is finitely
presented.
\end{lem}

\begin{proof}
Let $Q$ be the stabilizer in $L_{2}$ of the same vertex in $T_{2}$,
i.e. $V=P\cap Q$. Let $N$ be the (finite) fiber of $Q$, i.e. $\nicefrac{Q}{N}$
is a 2-orbifold group. The group $\nicefrac{V}{V\cap N}$ is a subgroup
of $\text{\ensuremath{\nicefrac{Q}{N}}}$, i.e. $\nicefrac{V}{V\cap N}$
is also a $2$-orbifold group, which is the orbifold group of the
cover of the original 2-orbifold. Since $\nicefrac{V}{V\cap N}$ is
finitely generated, the new 2-orbifold must be finite type, i.e. finitely
presented.

This means that $V$ is finitely presented subgroup of the $QH$ vertex
$Q$.
\end{proof}
\begin{lem}
Let $V$ be a virtually abelian vertex group in $T_{2}$, then $V$is
finitely generated virtually abelian, and thus finitely presented.
\end{lem}

\begin{proof}
Denote by $v$ the vertex in $T_{2}$ stabilized by $V$. Let $Q$
be the the stabilizer of $v$ in $L_{2}$ (as vertex in $T_{M}^{2}$).
The group $Q$ is virtually abelian (since $v$ is virtually abelian
vertex), and finitely generated relative to its edges stabilizer.
Since edge stabilizers are finitely generated so is $Q$. This means
that $V$ is a subgroup of a finitely generated virtually abelian
group so V itself is also finitely generated virtually abelian, and
thus finitely presented.
\end{proof}
This finishes the proof of \ref{prop:totally-rigid-step-1}.

Now that we have $Q_{1},\ldots,Q_{k}$ we can repeat the process with
each $Q_{i}$ separately, considering the quotient $L_{2}\onto L_{3}$,
and so on. This continue until we get to $L_{s}$. All the finitely
many subgroup we have obtained are contained in rigid subgroups of
$L_{s}$. In the case where the resolution is maximal (i.e. $L_{s}$
is non-divergent) we finish the process. In the case where $L_{s}$
is divergent we could apply the procedure one more time.

Wrapping this discussion we get
\begin{prop}
\label{prop:totally-rigid-subgroups}Let $G$ be finitely generated
group. Given a divergent limit group $L_{1}$ over $G$ which does
not admit an infinite resolution of proper limit quotients, then for
every resolution of shortening quotients 
\[
L_{1}\onto L_{2}\onto\cdots\onto L_{s}
\]
there exists finite set of finitely generated subgroups $P_{1},\ldots,P_{k}$
of $L_{1}$ such that
\begin{enumerate}
\item \label{enu:totally-rigid-lemma-1}$L_{1}$ is finitely presented relative
to $P_{1},\ldots,P_{k}$.
\item \label{enu:totally-rigid-lemma-2}The image of each $P_{1},\ldots,P_{k}$
in every group in the resolution (including $L_{1}$) is contained
in some rigid vertex group of the the JSJ decomposition of the corresponding
group, maybe aside from $L_{s}$ (which might be non-divergent, and
thus doesn't have a JSJ decomposition).
\item \label{enu:totally-rigid-lemma-3}Let $\left(\nu_{n}^{\left(i\right)}\right)_{n}$
be the defining sequence of $L_{i}$. For any finitely generated lift
$\tilde{P}_{d}$ of $P_{d}$ we have that for all $1\leq i,j\leq s$,
the maps $\restriction{\nu_{n}^{\left(i\right)}}{\tilde{P}_{d}}$
and $\restriction{\nu_{n}^{\left(j\right)}}{\tilde{P}_{d}}$ differ
by post-conjugation in $\Gamma$. 
\item \label{enu:totally-rigid-lemma-4}If $L_{s}$ is non-divergent (i.e.
when the resolution cannot be continued using shortening quotients)
then each of $P_{1},\ldots,P_{k}$ are a non-divergent limit group,
i.e. $P_{1},\ldots,P_{k}\in\lsbt\left(L_{i},\left(\nu_{n}^{\left(i\right)}\right)_{n}\right)$.
\end{enumerate}
\end{prop}

\begin{rem}
Notice that we could look at $P_{d}$ as limit group using any $\restriction{\nu_{n}^{\left(i\right)}}{\tilde{P}_{d}}$
as defining sequence (where $\tilde{P}_{d}\leq G$ is finitely generated
lift of $P_{d}$). Thanks to \enuref{totally-rigid-lemma-3}, it does
not matter which defining sequence we use to say that $P_{d}$ is
non-divergent.
\end{rem}

\begin{proof}
The proof of \enuref{totally-rigid-lemma-1} and \enuref{totally-rigid-lemma-2}
is the discussion gone through in the section. Item \enuref{totally-rigid-lemma-3}
follows from the fact that we obtained $L_{s}$ by shortening quotients,
thus one can apply \lemref{shortening_on_rigid_is_conj} iteratively.

Item \enuref{totally-rigid-lemma-4} follows from the fact that given
any finite generating set $\tilde{S}$ of $\tilde{P}_{d}$, the scaling
factor $\left\Vert \restriction{\nu_{n}^{\left(s\right)}}{\tilde{P}_{d}}\right\Vert _{\tilde{S}}$
can be bound by $K\left\Vert \nu_{n}^{\left(s\right)}\right\Vert _{S}$
where $S$ is a finite generating set of $G$ and $K=\max\left\{ \left|\tilde{s}\right|_{S}\mid\tilde{s}\in\tilde{S}\right\} $
where $\left|\cdot\right|_{S}$ is the word length in $G$ relative
to $S$. Now $\lim_{\omega}\left\Vert \restriction{\nu_{n}^{\left(s\right)}}{\tilde{P}_{d}}\right\Vert _{\tilde{S}}\leq\lim_{\omega}K\left\Vert \nu_{n}^{\left(s\right)}\right\Vert _{S}<\infty$.
\end{proof}

\subsection{\protect\label{subsec:Descending-Chains-Condition}Descending Chains
Condition}

The goal of this section is to prove that there is no sequence of
proper epimorphisms (satisfying some conditions) between limits groups.
This section is a generalization of \cite[Section 2]{JA-SE}, for
general hyperbolically acylindrical groups.
\begin{defn}
Let $G$ be finitely generated group. Define a relation on limit groups
which are defined over $G$ (i.e quotients of $G$) by $\left(L_{1},\varphi_{\infty}^{1}\right)\geq\left(L_{2},\varphi_{\infty}^{2}\right)$
if there is $\tau:L_{1}\onto L_{2}$ that commutes with $\varphi_{\infty}^{1}$
and $\varphi_{\infty}^{2}$ and that if $1\neq g\in B_{L_{1}}$ then
$1\neq\tau\left(g\right)\in B_{L_{2}}$ (see \defref{boundedness}).
\end{defn}

The following follows easily from the definition of shortening quotient
\begin{lem}
If $L_{1},\left(\varphi_{n}^{1}\right)$ is a limit group and $L_{2},\left(\varphi_{n}^{2}\right)$
is a shortening quotient of $L_{1}$ then $\left(L_{1},\varphi_{\infty}^{1}\right)\geq\left(L_{2},\varphi_{\infty}^{2}\right)$.
\end{lem}

\begin{proof}
In clear that $L_{1}\onto L_{2}$ commute with $\varphi_{\infty}^{1}$
and $\varphi_{\infty}^{2}$. Let $g\in B_{L_{1}}$ and $\tilde{g}\in G$
a $\varphi_{\infty}^{1}$ lift of $g$. Since $g$ is bounded trace
then $g$ is elliptic in the JSJ of $L_{1}$, thus from some point
onward, $\varphi_{n}^{1}\left(\tilde{g}\right)$ and $\varphi_{n}^{2}\left(\tilde{g}\right)$
differ by conjugation in $\Gamma$. This means that 
\begin{align*}
\lim_{\omega}\inf_{x\in X}d\left(x,\varphi_{n}^{2}\left(\tilde{g}\right)x\right) & =\lim_{\omega}\inf_{x\in X}d\left(x,c^{-1}\varphi_{n}^{1}\left(\tilde{g}\right)cx\right)=\\
 & =\lim_{\omega}\inf_{x\in X}d\left(cx,\varphi_{n}^{1}\left(\tilde{g}\right)cx\right)<\infty
\end{align*}
 thus $\varphi_{\infty}^{2}\left(\tilde{g}\right)\in B_{L_{2}}$.
\end{proof}
\begin{prop}
\label{prop:dcc}There exists no properly descending chain of limit
group (where proper means that the quotients in the chains are not
injective).
\end{prop}

\begin{proof}
Assume that there exists some infinite descending chain of limit groups
in this relation. Among all such infinite chains $\left(L_{1},\varphi_{\infty}^{1}\right)>\left(L_{2},\varphi_{\infty}^{2}\right)>\left(L_{3},\varphi_{\infty}^{3}\right)>\cdots$,
with corresponding defining sequences $\left(\varphi_{n}^{i}\right)$
(the sequence $\left(\varphi_{n}^{i}\right)$ define $L_{i}$, each
sequence runs over $n$ ) we choose the one satisfying the condition:
\begin{itemize}
\item[] If we have a different descending chain starting from $L_{i}$, 
\[
L_{i}>R_{i+1}>R_{i+2}>R_{i+3}>\cdots
\]
 then 
\[
\left|\ker\varphi_{\infty}^{i+1}\cap B_{i+1}\right|\geq\left|\ker\rho_{\infty}^{i+1}\cap B_{i+1}\right|
\]
where $\varphi_{\infty}^{i+1}$ it he limit quotient of $L_{i+1}$
and $\rho_{\infty}^{i+1}$ is the limit quotient of $R_{i+1}$ and
$B_{i+1}$ in the ball or radius $i+1$ around $1$ in $G$ (with
respect to a fixed finite generating set). This condition means that
among all infinite chains of limit groups we choose the one with the
``biggest kernels''.
\end{itemize}
Given a finite subset $S$ of $G$, denote $\left\Vert S\right\Vert _{n}^{i}=\inf_{x\in X}\max_{s\in S}d\left(x,\varphi_{n}^{i}\left(g\right)\cdot x\right)$.
Fixing $i$, we have that $\left<\varphi_{\infty}^{i}\left(S\right)\right>\in\sbt\left(L_{i}\right)$
if and only if $\lim_{w}\left\Vert S\right\Vert _{n}^{i}<\infty$
(see the definition of $B_{L_{i}}$, \ref{def:boundedness}). In case
$\left<\varphi_{\infty}^{i}\left(S\right)\right>\in B_{L_{i}}$ we
denote $\left\Vert S\right\Vert _{\infty}^{i}=\lim_{w}\left\Vert S\right\Vert _{n}^{i}$

\noindent Now for each $i$ we take $n_{i}$ such that 
\begin{itemize}
\item $\ker\varphi_{n_{i}}^{i}\cap B_{i}=\ker\varphi_{\infty}^{i}\cap B_{i}$
\item Given $S\subseteq B_{i}$ if $\left<\varphi_{\infty}^{i}\left(S\right)\right>\notin\sbt\left(L_{i}\right)$
then $\left\Vert S\right\Vert _{n_{i}}^{i}>i$.
\end{itemize}
i.e. we take $n_{i}$ such that $\varphi_{n_{i}}^{i}$ and $\varphi_{\infty}^{i}$
agree on which elements in the ball of radius $i$ are sent to $1$
and which don't, and likewise which finitely generated subgroup (with
generators in ball of radius $i$) are not bounded trace elements
(this is not precise since being bounded trace is not a property captured
be single homomorphism from $G$ to $\Gamma$, but the idea is such)
. From the definition of stable kernel and of bounded-trace subgroups
and the fact that $B_{i}$ is finite, there is such $n_{i}$.

We can choose the values of $n_{i}$ to be increasing.

Denote $\psi_{i}=\varphi_{n_{i}}^{i}$ and denote the resulting limit
group by $L_{\infty}$. By construction, $\ker\psi_{\infty}=\bigcup_{i=1}^{\infty}\ker\varphi_{\infty}^{i}$,
i.e. $L_{\infty}$ is the direct limit of all $L_{n}$'s.
\begin{lem}
For $g\in B_{L_{\infty}}$ and $\tilde{g}\in G$ such that $\psi_{\infty}\left(\tilde{g}\right)=g$,
there exists $i_{0}$ such that for $i>i_{0}$ we have $\varphi_{\infty}^{i}\left(\tilde{g}\right)\in B_{L_{i}}$.

Moreover $P\in\sbt\left(L_{\infty}\right)$, and finitely generated
$\psi_{\infty}$-lift $\tilde{P}\leq G$ then there exists $i_{0}$
such that for $i>i_{0}$ we have $\varphi_{\infty}^{i}\left(\tilde{P}\right)\subseteq B_{L_{i}}$.
\end{lem}

\begin{proof}
Denote $r=\left|\tilde{g}\right|$ the word length of $\tilde{g}$.
Since $g\in B_{L_{\infty}}$ we have $\left\Vert g\right\Vert _{\infty}^{\infty}\coloneqq\lim_{w}\left\Vert g\right\Vert _{i}^{\infty}<\infty$
where $\left\Vert g\right\Vert _{i}^{\infty}\coloneqq\inf_{x\in X}d\left(x,\psi_{i}\left(g\right)\cdot x\right)=\inf_{x\in X}d\left(x,\varphi_{n_{i}}^{i}\left(g\right)\cdot x\right)=\left\Vert \tilde{g}\right\Vert _{n_{i}}^{i}$.
Thus we have that $\wsurelu$ (in $i$) the inequality $\left\Vert \tilde{g}\right\Vert _{n_{i}}^{i}<\left\Vert g\right\Vert _{\infty}^{\infty}+1$
holds. For every $i>\left\Vert g\right\Vert _{\infty}^{\infty}+1,\left|\tilde{g}\right|$
(where $\left|\tilde{g}\right|$ is the word length of $\tilde{g}$
in $G$) we have $\wsurelu$ (in $i$) that
\[
\left\Vert \tilde{g}\right\Vert _{n_{i}}^{i}<\left\Vert g\right\Vert _{\infty}^{\infty}+1<i
\]
meaning that there exists $i_{0}$ such that $\varphi_{\infty}^{i_{0}}\left(\tilde{g}\right)$
must be in $B_{L_{i_{0}}}$ for some $i$. The condition regarding
bounded element in the definition of the relation implies that if
$\varphi_{\infty}^{i}\left(g\right)\in B_{L_{i}}$ then $\varphi_{\infty}^{i+1}\left(g\right)\in B_{L_{i+1}}$,
finishing the first part.

For $\left(2\right)$, Let $S$ be a generating set of $\tilde{P}$.
We define $\left\Vert S\right\Vert _{i}^{\infty}=\inf_{x\in X}\max_{s\in S}d\left(x,\psi_{i}\left(s\right)\cdot x\right)=\inf_{x\in X}\max_{s\in S}d\left(x,\varphi_{n_{i}}^{i}\left(s\right)\cdot x\right)=\left\Vert S\right\Vert _{n_{i}}^{i}$
and $\left\Vert S\right\Vert _{\infty}^{\infty}=\lim_{w}\left\Vert S\right\Vert _{i}^{\infty}$.
Since $P$ is in $\sbt\left(L_{\infty}\right)$ then $\left\Vert S\right\Vert _{\infty}^{\infty}<\infty$.
Meaning the inequality $\left\Vert S\right\Vert _{i}^{\infty}=\left\Vert S\right\Vert _{n_{i}}^{i}<\left\Vert S\right\Vert _{\infty}^{\infty}+1$
holds $\wsurelu$. The set of values of $i$ such that $i>\left\Vert S\right\Vert ^{\infty}+1,\max_{s\in S}\left|\tilde{s}\right|$
and $\left\Vert S\right\Vert _{n_{i}}^{i}<\left\Vert S\right\Vert _{\infty}^{\infty}+1$
is $\omega$ large (its an intersection of two $\omega$ large sets).
For each $i$ in that set we have $\left\Vert S\right\Vert _{n_{i}}^{i}<i$
and thus $\varphi_{\infty}^{i}\left(\tilde{P}\right)=\left<\varphi_{\infty}^{i}\left(S\right)\right>\in\sbt\left(L_{i}\right)$.
Take $i_{0}$ thus that $\varphi_{\infty}^{i_{0}}\left(\tilde{P}\right)\in\sbt\left(L_{i_{0}}\right)$.
Notice that $\varphi_{\infty}^{i_{0}+1}\left(\tilde{P}\right)=\tau_{i_{0}}\left(\varphi_{\infty}^{i_{0}}\left(\tilde{P}\right)\right)$
and since by the definition of the relation we have $\tau_{i_{0}}\left(B_{L_{i_{0}}}\right)\subseteq B_{L_{i_{0}+1}}$,
we have $\varphi_{\infty}^{i_{0}+1}\left(\tilde{P}\right)\subseteq B_{L_{i_{0}+1}}$
finishing the proof.
\end{proof}
Also note that $L_{\infty}$ admit only finite proper descending chains
of this relation. This is true since if there exists $L_{\infty}>R_{1}>R_{2}>\cdots$
an infinite chain with proper surjective maps then exists $g\in G$
such that $g$ maps to $1$ in $R_{1}$ but not in $L_{\infty}$.
If the length of $g$ is $n$ then $g$ wouldn't be sent to $1$ in
$L_{n+1}$ (since $L_{\infty}$ is the direct limit, and the way we
have build $L_{\infty}$) meaning the sequence $L_{n}\onto R_{1}\onto R_{2}\onto\cdots$
would contradict the ``maximality'' we required from $L_{1}\onto L_{2}\onto L_{3}\onto\cdots$
.

We now may apply \propref{totally-rigid-subgroups} on $L_{\infty}$
(notice that in \subsecref{Totally-rigid-subgroup} we said that $L$
does not admit any chain of surjection of limit groups, but since
\cite[Lemma 6.7]{GR-HU} used shortening quotients which are $>$-descending
we could apply \subsecref{Totally-rigid-subgroup} to $L_{\infty}$).

Let $L_{\infty}\onto R_{1}\onto R_{2}\onto\cdots$ be a sequence of
shortening quotients, then this sequence must stabilize. This means
the sequence has the form
\[
L_{\infty}\onto R_{1}\onto\cdots\onto R_{d}
\]
where $R_{d}$ is non-divergent limit group from \ref{lem:shortening-argument}.

Let $P_{1},\ldots,P_{k}$ be the totally rigid subgroups of $L_{\infty}$
relative to this finite chain (as in \propref{totally-rigid-subgroups}).
The modular automorphisms of the each limit group in the resolution
act by conjugation on totally rigid subgroups (since they are contained
in a rigid vertex group) hence totally rigid subgroups are also non-divergent
(with respect to the defining sequence of $L_{\infty}$), and thus
they are all contained in $B_{L_{\infty}}$.

Let $P$ be a totally rigid subgroup of $L_{\infty}$, and let $\tilde{P}$
be a finitely generated lift of $P$ to $G$. \propref{totally-rigid-subgroups}
tells us that $P\in\sbt\left(L_{\infty}\right)$, thus there is $i_{0}$
such that for $i>i_{0}$ we have $\varphi_{\infty}^{i}\left(\tilde{P}\right)\subseteq B_{L_{i}}$.
The quotients $\tau_{i}:L_{i}\onto L_{i+1}$ of the sequence are injective
restricted to $\tilde{P}$ , thus the image in $L_{n}$ must be isomorphic
to $P$ from some index onward (since $L_{\infty}$ is the direct
product of the $L_{i}$s then $P$ is the direct limit of $\varphi_{\infty}^{i}\left(\tilde{P}\right)s$).

The group $L_{\infty}$ is finitely presented relative to $P_{1},\ldots,P_{k}$,
and since $L_{\infty}$ is the direct limit of the $L_{s}$'s, there
exists an index from which onward all the other relations on $L_{\infty}$
hold.

This means that from some index all the other relations of $L_{\infty}$
hold and all the totally rigid limit group limit groups are embedded.
Hence from index onward point all the $L_{n}$ are quotients of $L_{\infty}$,
which can't happen, since $L_{\infty}$ is the direct limit of a sequence
of non-invective morphisms.
\end{proof}
As noted in the beginning of \subsecref{Totally-rigid-subgroup},
if a limit group admits no infinite sequence of descending limit groups
then all non bounded-trace subgroups are finitely generated, and together
with \propref{dcc} this shows the following lemma and its corollary
\begin{lem}
If $L$ is a divergent limit group, then all non bounded-trace virtually
abelian subgroups of $L$ are finitely generated.
\end{lem}

\begin{cor}
All edges in the modular decomposition of $L$ (and thus in its JSJ
decomposition) are finitely generated. Thus all vertex group in the
JSJ are finitely generated.
\end{cor}

\section{\protect\label{sec:Case-of-Product}Product Diagram Case}

\global\long\def\im{\text{Im}}%
\global\long\def\conj{\text{Conj}}%
\global\long\def\res{\text{Res}}%
\global\long\def\acts{\curvearrowright}%
\global\long\def\jsj{\text{JSJ}}%
\global\long\def\surject{\twoheadrightarrow}%
\global\long\def\comp{\text{Comp}}%
\global\long\def\sker{\underrightarrow{\ker}^{\omega}}%
\global\long\def\amal#1#2#3{#1\underset{#3}{*}#2}%
 
\global\long\def\hom{\text{Hom}}%
\global\long\def\isom{\text{Isom}}%
\global\long\def\bndv{\text{qTriv}}%
\global\long\def\bnd{\bndv_{c}}%
\global\long\def\walmost{\ \omega\text{-almost}}%
\global\long\def\wsurelu{\ \omega\text{-almost surely}}%
\global\long\def\onto{\twoheadrightarrow}%
\global\long\def\restriction#1#2{\left.\kern-\nulldelimiterspace#1\vphantom{\big|}\right|_{#2}}%
\global\long\def\mod{\text{Mod}}%
\global\long\def\aut{\text{Aut}}%
\global\long\def\stab{\text{Stab}}%
\global\long\def\by{\mathbf{Y}}%
 
\global\long\def\projcopx{\mathcal{P}_{K}}%
\global\long\def\diam{\text{diam}}%
\global\long\def\qtms{\mathcal{C}_{K}}%
\global\long\def\c{\mathcal{C}}%
\global\long\def\cal#1{\mathcal{#1}}%

\subsection{Basic Definitions}

A product diagram is a space of the form $\mathcal{X}=V_{1}\times\cdots\times V_{m}$
where each of the spaces is an unbounded hyperbolic space. We equip
$\mathcal{X}$ with the $l_{1}$ metric.
\begin{defn}
Let $\sigma\in S_{m}$ be a permutation and $f_{i}:V_{i}\to V_{\sigma\left(i\right)}$
be an isometry for each $i\in\left\{ 1,\ldots,m\right\} $. We construct
an isometry of $\mathcal{X}$ by 
\[
\left(x_{1},\ldots,x_{m}\right)\mapsto\left(f_{\sigma^{-1}\left(1\right)}\left(x_{\sigma^{-1}\left(1\right)}\right),\dots,f_{\sigma^{-1}\left(m\right)}\left(x_{\sigma^{-1}\left(m\right)}\right)\right)
\]
An automorphism of (the product diagram) $\mathcal{X}$ is an isometry
obtained by such construction.The group $\aut\left(\mathcal{X}\right)$
is defined as all isometries constructed this way. Notice that there
is a natural morphism $\aut\left(\mathcal{X}\right)\to S_{m}$ which
extracts the underlying permutation.

An action of a group $\Gamma$ on $\mathcal{X}$ is a morphism $\Gamma\to\aut\left(\mathcal{X}\right)$
(note that this is more restrictive than action by isometries). We
define the \textbf{stabilizer of $\mathbf{V_{i}}$} by
\[
\stab\left(V_{i}\right)=\left\{ g\in\Gamma\mid\text{the associated permutation of \ensuremath{g} fixes \ensuremath{i}}\right\} 
\]
i.e. $\stab\left(V_{i}\right)$ is the stabilizer of $i$ in the induced
action of $\Gamma$ on $\left\{ 1.,\ldots,m\right\} $. 
\end{defn}

\begin{rem*}
This definition is analogous to the definition of automorphisms of
a HHS, and in fact a special case of it.
\end{rem*}
\begin{defn}
\defref{strictly-acylindrical-product}\label{def:strictly-acylindrical-product}The
action of a finitely generated group $\Gamma$ on product diagram
$\mathcal{X}$ is called \textbf{strictly acylindrical} if for every
$i=1,\ldots,m$ the action of $\stab\left(V_{i}\right)$ on $V_{i}$
is acylindrical.
\end{defn}

Fix a f.g. group $\Gamma$. Assume that that $\Gamma$ acts geometrically
and strictly acylindrically on $\mathcal{X}$. Let $\Gamma'\leq\Gamma$
be the kernel of the induced action of $\Gamma$ on $\left\{ 1,\ldots,m\right\} $,
i.e. $\Gamma'=\bigcap_{i=1}^{m}\stab\left(V_{i}\right)$. Notice that
$\Gamma'$ has finite index in $\Gamma$, and that $\Gamma'$ also
acts geometrically on $\mathcal{X}$. 

In the proof of \Thmref{product-diagram-noetherian} we will see that
it does not restrict the generality of the result to assume that $\Gamma=\Gamma'$.
In doing so we actually assume that $\stab\left(V_{i}\right)=\Gamma$
for all $i$.

Given a sequence $\left(\varphi_{n}\right)\in\hom\left(G,\Gamma\right)$,
we can look at the sequence in a few ways. The first is by looking
at the sequence of actions of $G$ on $\mathcal{X}$ induced by $\left(\varphi_{n}\right)$.
The others are by looking at the sequence of actions of $G$ on each
$V_{i}$ induced by $\left(\varphi_{n}\right)$, i.e. precompose with
$\left(\varphi_{n}\right)$ and then restrict the action to $V_{i}$.

From each induced series of actions we get its scaling factor, which
we denote $\left\Vert \varphi_{n}\right\Vert _{\mathcal{X}},\left\Vert \varphi_{n}\right\Vert _{1},\ldots,\left\Vert \varphi_{n}\right\Vert _{m}$
respectively.

We say that $\left(\varphi_{n}\right)$ is divergent/non-divergent
relative to $\mathcal{X}$ (resp. to $V_{i}$) if the induced sequence
of action on $\mathcal{X}$ (resp. $V_{i}$) is divergent/non-divergent
as in \ref{def:divergent_seq}. In short, we will say that $\left(\varphi_{n}\right)$
is $j$-divergent/$j$-non-divergent for divergent/non-divergent relative
to $V_{j}$.
\begin{lem}
\label{lem:non-divergent-rel-is-non-divergent}The sequence $\left(\varphi_{n}\right)$
is non-divergent relative to $\mathcal{X}$ if and only if it is $j$-non-divergent
for all $j$, $1\leq j\leq m$.
\end{lem}

\begin{proof}
Assume that $\left(\varphi_{n}\right)$ is non divergent relative
to $\mathcal{X}$, i.e. there exists some $D>0$ such that $\lim_{\omega}\left\Vert \varphi_{i}\right\Vert _{\mathcal{X}}<D$.
Namely we can take $x^{\left(n\right)}=\left(x_{1}^{\left(n\right)},\ldots,x_{m}^{\left(n\right)}\right)\in\mathcal{X}$
such that $d_{\mathcal{X}}\left(\varphi_{n}\left(s\right)\cdot x^{\left(n\right)},x^{\left(n\right)}\right)<D$
$\wsurelu$. Each sequence $\left(x_{i}^{\left(n\right)}\right)$
shows us that $\varphi_{n}$ is non-divergent relative to $V_{i}$.
We have $\wsurelu$ for all $s$ in the pre-chosen finite generating
set of $G$ and every $i$ in $\left\{ 1,\ldots,m\right\} $ that
\begin{align*}
d_{V_{i}}\left(\varphi_{n}\left(s\right)\cdot x_{i}^{\left(n\right)},x_{i}^{\left(n\right)}\right)\leq & \sum_{i=1}^{m}d_{V_{i}}\left(\varphi_{n}\left(s\right)\cdot x_{i}^{\left(n\right)},x_{i}^{\left(n\right)}\right)=\\
= & d_{\mathcal{X}}\left(\varphi_{n}\left(s\right)\cdot x^{\left(n\right)},x^{\left(n\right)}\right)<D
\end{align*}
thus $\wsurelu$ we have $\max\limits_{s\in S}d\left(\varphi_{n}\left(s\right)\cdot x_{i}^{\left(n\right)},x_{i}^{\left(n\right)}\right)<D$.
So $\lim_{\omega}\left\Vert \varphi_{n}\right\Vert _{i}<D$ for all
$i$ in $\left\{ 1,\ldots,m\right\} $

Assume that $\left(\varphi_{n}\right)$ is non divergent relative
to $V_{1},\ldots,V_{m}$.

Let $D>0$ such that $\lim_{\omega}\left\Vert \varphi_{i}\right\Vert _{i}<D$
for all $i$. For each $i$, let $\left(x_{i}^{\left(n\right)}\right)\in V_{i}$
be a sequence such that $\wsurelu$ for all $s$ in the generating
set $d_{V_{i}}\left(\varphi_{n}\left(s\right)\cdot x_{i}^{\left(n\right)},x_{i}^{\left(n\right)}\right)<D$.
The point $x^{\left(n\right)}=\left(x_{1}^{\left(n\right)},\ldots,x_{m}^{\left(n\right)}\right)\in\mathcal{X}$
shows us that $\varphi_{n}$ is non-divergent relative to $\mathcal{X}$.
\[
d_{\mathcal{X}}\left(\varphi_{n}\left(s\right)\cdot x^{\left(n\right)},x^{\left(n\right)}\right)=\sum_{i=1}^{m}d_{V_{i}}\left(\varphi_{n}\left(s\right)\cdot x_{i}^{\left(n\right)},x_{i}^{\left(n\right)}\right)<mD
\]
and this inequality is true $\wsurelu$ for all $s$ in the generating
set.
\end{proof}

\subsection{Higher rank Resolution over a Product\protect\label{subsec:Higher-rank-Resolution-product}}

Let $\left(\varphi_{n}\right)\subset\hom\left(G,\Gamma\right)$ be
a divergent sequence (relative to $\mathcal{X}$).

Let $V_{j}$ be a factor such that $\left(\varphi_{n}\right)$ is
$j$-divergent. Regarding only the action of $\Gamma$ on $V_{j}$,
one can use \cite{GR-HU} to construct a resolution of shortening
quotients. i.e. let $L_{1}^{j}=\nicefrac{G}{\sker\varphi_{n}}$, and
$L_{2}^{j}$ a shortening quotient of $L_{1}^{j}$ (as described in
\subsecref{Modular-automorphism}) and so on. As in \propref{dcc},
this procedure must stabilize, i.e we get a sequence of shortening
quotients

\[
L_{1}^{j}\onto L_{2}^{j}\onto\cdots\onto L_{t_{j}}^{j}
\]
where $L_{t_{j}}^{j}$ is a $j$-non-divergent limit group (\lemref{shortening-argument}).
We call this sequence the $j$-resolution of $\left(\varphi_{n}\right)$.
For each limit group in this resolution there is the associated JSJ
decomposition $\Delta_{k}^{j}$, where $1\leq k\leq t_{j}$ (as discussed
in \subsecref{JSJ-decomposition}).

For each $j$ such that $\left(\varphi_{n}\right)$ is $j$-divergent
we create the corresponding resolution $L_{1}^{j}\onto L_{2}^{j}\onto\cdots\onto L_{t_{j}}^{j}$.
For values of $j$ where $\left(\varphi_{n}\right)$ is $j$-non-divergent
we set the resolution to be single group $L_{1}^{j}=\nicefrac{G}{\sker\varphi_{n}}$,
and the ``JSJ'' $\Delta_{1}^{j}$ to be the trivial splitting (i.e
$\Delta_{1}^{j}$ is a point). For those values of $j$, we say that
the resolution is trivial/single-level.

Note that even though independently of the value of $j$ the groups
$L_{1}^{j}$ are all the same (algebraically they are just $\nicefrac{G}{\sker\varphi_{n}}$)
and seem to be independent of the action of $V_{j}$, the JSJ decompositions,
and all other limit groups in the resolution do depend on the action
of $\Gamma$ on each $V_{j}$.

Conceptually, $\nicefrac{G}{\sker\varphi_{n}}$ and $L_{1}^{j}$ play
different roles. In order to prove noetherianity of $\Gamma$ we want
$\left(\varphi_{n}\right)$ to $\wsurelu$ factor through $\nicefrac{G}{\sker\varphi_{n}}$,
it plays an ``algebraic role'', but $L_{1}^{j}$ plays a ``geometric
role''. When one tries to relax the strict acylindricity assumptions
to weak acylindricity (see \cite{Se23}), the action of the limit
groups on the limit real tree might have ``big'' kernels, thus one
need to change the definition of limit groups, making the isomorphism
between $L_{1}^{j}$ and $\nicefrac{G}{\sker\varphi_{n}}$ no longer
true. We heavily use this isomorphism in the proof of \Thmref{product-diagram-noetherian},
using the fact that the JSJ of each $L_{1}^{j}$ is also a decomposition
of $\nicefrac{G}{\sker\varphi_{n}}$.

\begin{figure}[h]
\[
\xymatrix{ &  & G\ar@{->>}[d]\\
 &  & \nicefrac{G}{\sker\left(\varphi_{n}\right)}\ar@{->>}[lld]\ar@{->>}[ld]\ar@{->>}[rd]\ar@{->>}[rrd]\\
L_{1}^{1}\ar@{->>}[d] & L_{1}^{2}\ar@{->>}[d] &  & L_{1}^{m-1}\ar@{->>}[d] & L_{1}^{m}\ar@{->>}[d]\\
L_{2}^{1}\ar@{->>}[d] & L_{2}^{2}\ar@{->>}[d] &  & L_{2}^{m-1}\ar@{->>}[d] & L_{2}^{m}\ar@{->>}[d]\\
\vdots\ar@{->>}[d] & \vdots\ar@{->>}[d] &  & \vdots\ar@{->>}[d] & \vdots\ar@{->>}[d]\\
L_{t_{1}}^{1} & L_{t_{2}}^{2} &  & L_{t_{m-1}}^{m-1} & L_{t_{m}}^{m}\\
\mbox{{\tiny resolution relative to \ensuremath{V_{1}}}} & \text{{\tiny resolution relative to \ensuremath{V_{2}}}} &  & \text{{\tiny resolution relative to \ensuremath{V_{m-1}}}} & \text{{\tiny resolution relative to \ensuremath{V_{m}}}}
}
\]

\caption{\protect\label{fig:higher-rank-resolution}The higher rank resolution
associated to $\left(\varphi_{n}\right)$}
\end{figure}
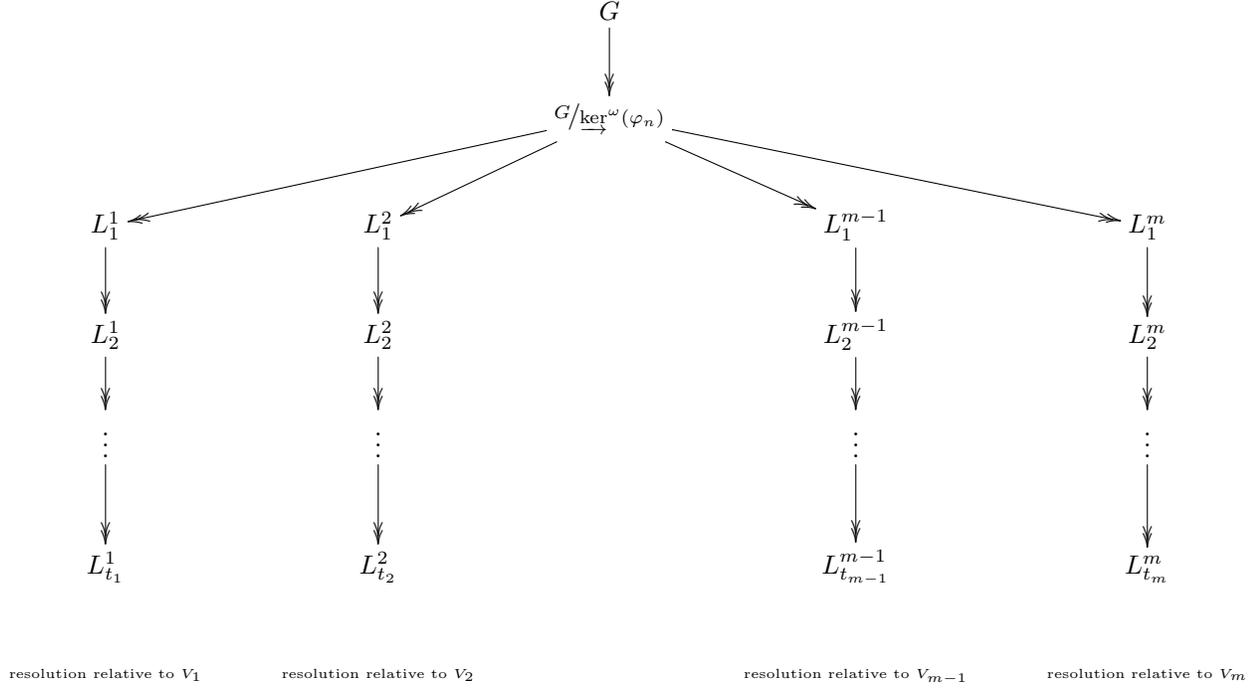

We now arrange all the different resolutions (some of which might
be trivial) in a diagram as in \figref{higher-rank-resolution} (notice
we distinguish between the first level of each resolution and $\nicefrac{G}{\sker\varphi_{n}}$).

\begin{defn}
\label{def:higer-rank-res}We call this diagram the \textbf{higher
rank resolution of $\bm{\left(\varphi_{n}\right)}$}.
\end{defn}

\subsection{Proof of Noetherianity}

We will prove that every sequence $\left(\varphi_{n}\right)\in\hom\left(G,\Gamma\right)$
factors through its corresponding limit groups by induction on the
size of the set $\left\{ j\mid\varphi_{n}\text{ is divergent relative to \ensuremath{V_{j}}}\right\} $.
If the size of this set is $k$, then we say that $\left(\varphi_{n}\right)$
has \textbf{$k$ divergent factors}. If $\left(\varphi_{n}\right)$
is non divergent relative to $\mathcal{X}$ then it has $0$ divergent
factors.
\begin{lem}[Induction Base]
\label{lem:base-induction-product}All non-divergent (relative to
$\mathcal{X}$) sequences $\wsurelu$ factor through their corresponding
limit group.
\end{lem}

\begin{proof}
Let $\left(\varphi_{n}\right)\in\hom\left(G,\Gamma\right)$ be such
sequence, i.e. a sequence such that $\lim_{\omega}\left\Vert \varphi_{n}\right\Vert <D$
for some $D>0$.

Let $x_{n}$ be a minimally displaced point of $\varphi_{n}$, meaning
$x_{n}\in\mathcal{X}$ such that $\max_{s\in S}d_{\mathcal{X}}\left(x_{n},\varphi_{n}\left(s\right)x_{n}\right)<D$
for each $n$. Since the action on $\mathcal{X}$ is co-compact, once
can assume, up to post-compose $\varphi_{n}$ with conjugation, that
all of the $x_{n}$ are in same compact set. Choose $x_{0}$ in this
compact set. Since $x_{0}$ has uniform bound on the distances form
every $x_{n}$, it also has bounded displacement distance. i.e.
\[
\max_{s\in S}d_{{\cal X}}\left(x_{0},\varphi_{n}\left(s\right)x_{0}\right)<D+K
\]
 where $K$ depend on the radius of the compact set containing all
$x_{n}$. 

This means that there exists an $\omega$ large set such that for
all $s\in S$ (the fixed finite generating set of $G$) and all $n$
in the $\omega$ large set we have that $\varphi_{n}\left(s\right)x_{0}$
is in a $\left(D+K\right)$-ball around $x_{0}$. Since that action
in properly discontinuous then $\varphi_{n}\left(s\right)$ can only
lie in a finite set, meaning that $\left(\varphi_{n}\right)$, up
to post-conjugation, consist of only finitely many morphisms (since
$S$ generates $G$). This means that there exists an $\omega$ large
set such that the maps $\left\{ \varphi_{n}\right\} $ are all conjugate.
Restricting to this $\omega$-large set, $\ker\varphi_{n}$ is constant
(thus equal to $\sker\varphi_{n}$), meaning every $\varphi_{n}$
(for $n$ in this $\omega$-large set) factors through $\nicefrac{G}{\sker\varphi_{n}}$.

\end{proof}
\begin{prop}[Induction step]
\label{prop:step-induction-product}Assume that for some $k\leq m$,
every sequence $\left(\varphi_{n}\right)\in\hom\left(G,\Gamma\right)$
with at most $k$ divergent factors $\wsurelu$ factors through its
corresponding limit group.

Then every sequence $\left(\varphi_{n}\right)\in\hom\left(G,\Gamma\right)$
with at most $k+1$ divergent factors $\wsurelu$ factors through
its corresponding limit group.
\end{prop}

\begin{proof}
Let $\left(\varphi_{n}\right)\in\hom\left(G,\Gamma\right)$ be a sequence
with $k+1$ divergent factors. We can rearrange the factor such that
for $1\leq j\leq k+1$ the sequence $\left(\varphi_{n}\right)$ is
$j$-divergent and for $k+1<j$ the sequence $\left(\varphi_{n}\right)$
is $j$-non-divergent.

Let
\[
L_{1}^{1}\to\cdots\to L_{t_{1}}^{1}
\]
be the resolution associated with $V_{1}$ (as in \ref{fig:higher-rank-resolution}).
Denote by $\left(\nu_{n}\right)\in\hom\left(G,\Gamma\right)$ the
defining sequence of $L_{t_{1}}^{1}$.

Let $P_{1},\ldots,P_{s}$ be the totally rigid subgroups of $L_{1}^{1}$
(see \propref{totally-rigid-subgroups}). Since each totally rigid
subgroup is finitely generated we can choose finitely generated $\varphi_{\infty}$-lifts
$\tilde{P}_{1},\ldots,\tilde{P}_{s}\leq G$ of $P_{1},\ldots,P_{s}$. 

Let $P$ be one of the $P_{1},\ldots,P_{s}$, and denote its finitely
generated lift by $\tilde{P}$.

\Propref{totally-rigid-subgroups} tells us that $\left(\restriction{\varphi_{n}}{\tilde{P}}\right)_{n}$
and $\left(\restriction{\nu_{n}}{\tilde{P}}\right)_{n}$ differ only
by conjugation and that $\left(\restriction{\varphi_{n}}{\tilde{P}}\right)_{n}$
is non-divergent relative to $V_{1}$ (since we took the resolution
such that $L_{t_{1}}^{1}$ is non-divergent relative to $V_{1}$).
Also, if $\left(\varphi_{n}\right)_{n}$ is non-divergent relative
to some $V_{j}$ then $\left(\restriction{\varphi_{n}}{\tilde{P}}\right)_{n}$
is also non-divergent relative to $V_{j}$.

All in all, $\left(\restriction{\varphi_{n}}{\tilde{P}}\right)_{n}:\tilde{P}\to\Gamma$
could be divergent relative to $V_{i}$ where $1<i\leq k+1$, i.e.
$\left(\restriction{\varphi_{n}}{\tilde{P}}\right)_{n}$ has at most
$k$ divergent components. Thus, by induction assumption, we get $\left(\restriction{\varphi_{n}}{\tilde{P}}\right)_{n}$
factor through $\restriction{\varphi_{\infty}}{\tilde{P}}$ for $\walmost$
any $n$.

This is true for all $\restriction{\varphi_{n}}{\tilde{P}_{1}},\ldots,\restriction{\varphi_{n}}{\tilde{P}_{s}}$
, i.e $\restriction{\varphi_{n}}{\tilde{P}_{j}}$ $\wsurelu$ factor
through $\restriction{\varphi_{\infty}}{\tilde{P}_{j}}$ for every
$1\leq j\leq s$ . Using \lemref{relative_presented_lemma} we get
that $\left(\varphi_{n}\right)_{n}$ $\walmost$ factor-through $L_{1}^{1}=\nicefrac{G}{\sker}\varphi_{n}$.
\end{proof}
Combining the two statements \Lemref{base-induction-product} and
\Propref{step-induction-product}, we conclude the main theorem:
\begin{thm}
\label{thm:product-diagram-noetherian}Any finitely generated groups
acting geometrically and strictly acylindrically on a product of hyperbolic
spaces is equationally noetherian
\end{thm}

\begin{proof}
Let $\hat{\Gamma}$ be the finite index subgroup of $\Gamma$ fixing
the factors. The group $\hat{\Gamma}$ still acts strictly acylindrically
and geometrically on the same product diagram. From \Lemref{base-induction-product}
and \Propref{step-induction-product} any sequence of morphism from
any finitely generated group $G$ to $\hat{\Gamma}$ $\wsurelu$ factor
through the corresponding limit group, thus $\hat{\Gamma}$ is equationally
noetherian, and virtually noetherian groups are noetherian groups
(\cite{BMR}), thus $\Gamma$ is noetherian.
\end{proof}

\section{\protect\label{sec:general-digram}General Diagram Case}

\global\long\def\im{\text{Im}}%
\global\long\def\conj{\text{Conj}}%
\global\long\def\res{\text{Res}}%
\global\long\def\acts{\curvearrowright}%
\global\long\def\jsj{\text{JSJ}}%
\global\long\def\surject{\twoheadrightarrow}%
\global\long\def\comp{\text{Comp}}%
\global\long\def\sker{\underrightarrow{\ker}^{\omega}}%
\global\long\def\amal#1#2#3{#1\underset{#3}{*}#2}%
 
\global\long\def\hom{\text{Hom}}%
\global\long\def\isom{\text{Isom}}%
\global\long\def\bndv{\text{qTriv}}%
\global\long\def\bnd{\bndv_{c}}%
\global\long\def\wsurelu{\ \omega\text{-almost surely}}%
\global\long\def\onto{\twoheadrightarrow}%
\global\long\def\restriction#1#2{\left.\kern-\nulldelimiterspace#1\vphantom{\big|}\right|_{#2}}%
\global\long\def\mod{\text{Mod}}%
\global\long\def\aut{\text{Aut}}%
\global\long\def\stab{\text{Stab}}%
\global\long\def\by{\mathbf{Y}}%
\global\long\def\into{\hookrightarrow}%
\global\long\def\out{\text{Out}}%
 
\global\long\def\projcopx{\mathcal{P}_{K}}%
\global\long\def\diam{\text{diam}}%
\global\long\def\qtms{\mathcal{C}_{K}}%
\global\long\def\c{\mathcal{C}}%
\global\long\def\sfrak{\mathfrak{S}}%

In this section $\Gamma$ will be an HHG. Our main result requires
some extra assumption on $\Gamma$ (strict acylindricity and colorability)
that we will define soon.

The proof of noetherianity in the general case is essentially the
same as in the case of a product of hyperbolic groups, but it requires
additional tools.

The main objective of the first two subsections is to explain how
one can see colorable HHGs (see \defref{colorable-HHG}) as acting
on a product of hyperbolic spaces. \Subsecref{Projection-Complexes}
is a quick survey of the projection complex machinery defined by Besvina,
Bromberg, Fujiwara and Sisto \cite{BBF,BBFS}. In \subsecref{Colorable-HHGs}
we use the machinery to view colorable HHGs as acting on product of
hyperbolic spaces.

\subsection{\protect\label{subsec:Projection-Complexes}Projection Complexes,
quasi-trees of metric spaces}
\begin{defn}[{\cite[Section 2]{BBFS}}]
Let $\by$ be a set such that for each $Y\in\by$ there is a metric
space $\left(\mathcal{C}Y,\rho_{Y}\right)$ and for every $Y\in\by$
a function $\pi_{Y}$ from $\by\backslash\left\{ Y\right\} $ to subsets
of $\mathcal{C}Y$, called \textbf{projection}. For every $Y\in\by$
and $X,Z\in\by\backslash\left\{ Y\right\} $ we define the corresponding
pseudo-distance function
\[
d_{Y}^{\pi}\left(X,Z\right)=\diam\left(\pi_{Y}\left(X\right)\cup\pi_{Y}\left(Z\right)\right)
\]

When $\left\{ \pi_{Y}\right\} $ is clear we write $d_{Y}$ for $d_{Y}^{\pi}$.

We call the set of functions $\left\{ d_{Y}\right\} $ a \textbf{strong
projection} if for some $\theta\geq0$ it satisfy the \textbf{\textit{strong
projection axioms}}:
\begin{enumerate}
\item[\textbf{(SP 1)}] $d_{Y}\left(X,Z\right)=d_{Y}\left(Z,X\right)$
\item[\textbf{(SP 2)}] $d_{Y}\left(X,Z\right)+d_{Y}\left(Z,W\right)\geq d_{Y}\left(X,W\right)$
\item[\textbf{(SP 3)}] If $d_{Y}\left(X,Z\right)>\theta$ then $d_{Z}\left(X,W\right)=d_{Z}\left(Y,W\right)$
for all $W\in\mathbf{Y}\backslash\left\{ Z\right\} $
\item[\textbf{(SP 4)}] $d_{Y}\left(X,X\right)\leq\theta$
\item[\textbf{(SP 5)}] $\#\left\{ Y\mid d_{Y}\left(X,Z\right)>\theta\right\} $ is finite
for all $X,Z\in\mathbf{Y}$
\end{enumerate}
\end{defn}

\begin{defn}[{\cite[Lemma 2.2]{BBFS}}]
Let $\by$, $\left\{ {\cal C}Y\right\} $, $\left\{ \pi_{Y}\right\} $
and $\left\{ d_{Y}\right\} $ be a strong projections, an let $K>2\theta$.
Define $Y_{K}\left(X,Z\right)=\left\{ Y\in\mathbf{Y}\mid d_{Y}\left(X,Z\right)>K\right\} $
and a partial order on $Y_{K}\left(X,Z\right)$ by $Y_{0}<Y_{1}$
if the following equivalent conditions (\cite[Lemma 2.2]{BBFS}) hold:
\begin{enumerate}
\item $d_{Y_{0}}\left(X,Y_{1}\right)>\theta$;
\item $d_{Y_{1}}\left(Y_{0},W\right)=d_{Y_{1}}\left(X,W\right)$ for all
$W\neq Y_{1}$;
\item $d_{Y_{1}}\left(X,Y_{0}\right)\leq\theta$;
\item $d_{Y_{1}}\left(Y_{0},Z\right)>\theta$;
\item $d_{Y_{0}}\left(W,Y_{1}\right)=d_{Y_{0}}\left(W,Z\right)$ for all
$W\neq Y_{0}$
\item $d_{Y_{0}}\left(Y_{1},Z\right)\leq\theta$
\end{enumerate}
\noindent We define the the \textbf{projection complex $\projcopx\left(\by\right)$
}to be the graph defined with the vertex set as $\mathbf{Y}$ and
$X,Z\in\mathbf{Y}$ are connected if $Y_{K}\left(X,Z\right)=\emptyset$.

\end{defn}

\begin{fact}
If $G$ is a group acting on $\by$ preserving $\left\{ \pi_{Y}\right\} $
(meaning $\pi_{gY}\left(gX\right)=\pi_{Y}\left(X\right)$) then the
construction of $\projcopx\left(\by\right)$ is $G$-equivariant,
and $G$ acts by isometries on $\projcopx\left(\by\right)$.
\end{fact}

\begin{thm}[{\cite[Theorem 3.5]{BBFS}}]
\label{thm:projection_complex_quasi_tree}For $K>3\theta$, $\projcopx\left(\by\right)$
is a quasi-tree.
\end{thm}

\begin{rem}
In \cite{BBF,BBFS} the authors build the projection complex without
the projection $\left\{ \pi_{Y}\right\} $. We choose to include them
in the definition of strong projection for simplicity, as they heavily
used in the sequel and we do not need the extra level of abstraction.
\end{rem}

\begin{defn}
Let $\by,\left\{ \mathcal{C}Y\right\} ,\left\{ \pi_{Y}\right\} $
and $\left\{ d_{Y}\right\} $ satisfy the strong projection axioms,
and let $K>0$. We define the \textbf{quasi-tree of metric spaces}
$\mathcal{C}_{K}\left(\by\right)$ (not to be confused with $\mathcal{C}Y$)
to be the disjoint union of the $\mathcal{C}Y$ when if $d_{\projcopx\left(\by\right)}\left(X,Y\right)=1$
we connect every point of $\pi_{X}\left(Z\right)$ to every point
of $\pi_{Z}\left(X\right)$ with a segment of length $K$.

As a convention, points in $\qtms\left(\by\right)$ will be denoted
by small letter, and the corresponding spaces in $\by$ containing
them would be denoted by the corresponding upper case letters (e.g.
$x,z\in\qtms\left(\by\right)$ and $X,Z\in\by$ indicate that $x\in{\cal C}X$
and $x\in{\cal C}Z$).
\end{defn}

\begin{rem}
One can think of a quasi-tree of metric spaces as a blow up of $\projcopx\left(\by\right)$,
where we blow every $Y$ to $\mathcal{C}Y$ and connect it to its
neighbors at the images of projections, see \Figref{Illustration-of-quasi-tree}. 

\begin{figure}
\centering{}\includegraphics[scale=0.2]{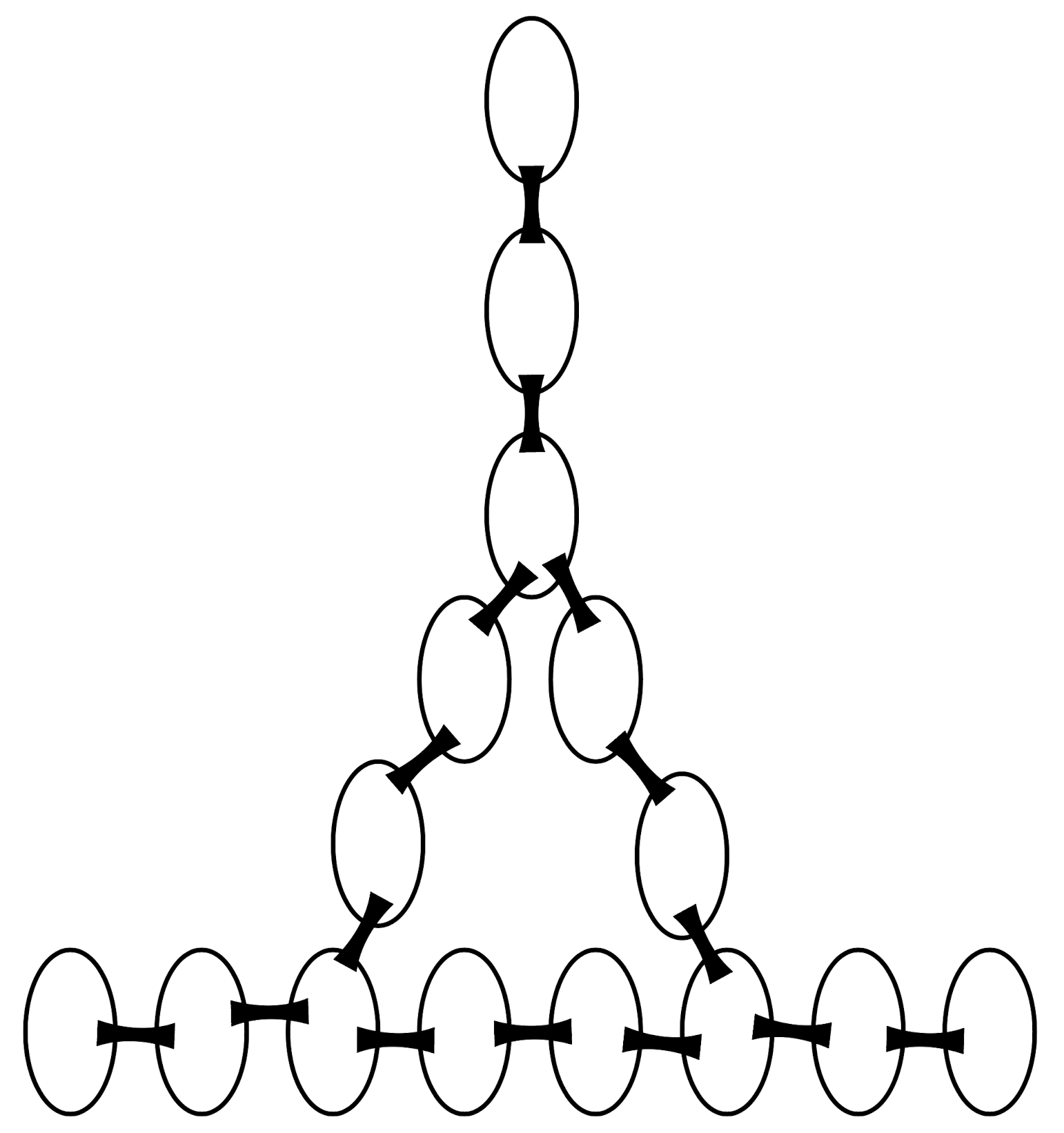}\caption{\protect\label{fig:Illustration-of-quasi-tree}Illustration of quasi-tree
of metric spaces as blow up of $\protect\projcopx\left(\protect\by\right)$}
\end{figure}

In light of the the last remark, we have similar result to \ref{thm:projection_complex_quasi_tree},
which state that quasi-tree of metric spaces inherit hyperbolic properties
of the spaces used to construct it.
\end{rem}

\begin{thm}[{\cite[Theorem 4.17]{BBF}}]
Let $\by$, $\left\{ {\cal C}Y\right\} $, $\left\{ \pi_{Y}\right\} $
and $\left\{ d_{Y}\right\} $ be strong projections such that each
$\c Y$ is $\delta$-hyperbolic then $\qtms\left(\by\right)$ is hyperbolic.
\end{thm}

\begin{defn}
Let $\pi_{Y}^{\flat}:\mathcal{C}_{K}\left(\by\right)\to\c Y$ be the
map defined by
\[
\pi_{Y}^{\flat}\left(x\right)=\begin{cases}
x & X=Y\\
\pi_{Y}\left(X\right) & X\neq Y
\end{cases}
\]
where $X,Y\in\by$ and $X$ is the space containing $x$. This map
coarsely equal to the closest point projection of $\mathcal{C}_{K}\left(\by\right)$
onto $\c Y$ (\cite[4.10]{BBF}), thus we call $\pi_{Y}^{\flat}$
the \textbf{closest point projection}. If a group $G$ acts on $\by$
then $\pi_{Y}$ is $G$-equivariant.

For $x,z\in\qtms\left(\by\right)$ and $Y\in\by$ we define
\[
d_{Y}\left(x,z\right)=\diam\left(\pi_{Y}^{\flat}\left(x\right)\cup\pi_{Y}^{\flat}\left(z\right)\right)
\]
and $Y_{K}\left(x,z\right)=\left\{ Y\in\by\mid d_{Y}\left(x,z\right)>K\right\} $.
\end{defn}

\begin{rem}
Notice that the definition of $Y_{K}\left(x,z\right)$ and $Y_{K}\left(X,Z\right)$
(for $X,Z\in\by$) differ since $Y_{K}\left(x,z\right)$ might contain
$X$ if $x$ and $\pi_{X}^{\flat}\left(z\right)$ are far in ${\cal C}X$
(respectively for $Z$). This can not happen in $Y_{K}\left(X,Z\right)$.
\end{rem}

We have a distance estimation formula in $\mathcal{C}_{K}\left(\by\right)$
using the projections of points in $\mathcal{C}_{K}\left(\by\right)$. 
\begin{thm}[{\cite[6.3]{BBFS}}]
Let $K$ be sufficiently large. Then given $x,z\in\mathcal{C}_{K}\left(\by\right)$
where $x\in X$ and $z\in Z$ ($X,Z\in\by$) we have
\[
\frac{1}{4}\sum_{Y\in Y_{K}\left(x,z\right)}d_{Y}\left(x,z\right)\leq d_{\mathcal{C}_{K}\left(\by\right)}\left(x,z\right)\leq2\sum_{Y\in Y_{K}\left(x,z\right)}d_{Y}\left(x,z\right)+3K
\]
\end{thm}

\begin{rem}
This show that the inclusion $\c Y\into\mathcal{C}_{K}\left(\by\right)$
is quasi-isometry embedding.
\end{rem}

\subsection{\protect\label{subsec:Colorable-HHGs}Colorable HHGs}

Let $\mathcal{X}$ be a HHS, as defined in \cite{BHS1,BHS2,S}. Let
$\mathfrak{S}$ be its index/domain set.

From any subset $\mathfrak{L}\subset\mathfrak{S}$ of pairwise transverse
elements, and $X,Y\in\mathfrak{L}$ there is a projection set $\rho_{Y}^{X}\subset Y$
uniformly bounded. If we define $\pi_{Y}\left(X\right)=\rho_{Y}^{X}$
we get a pseudo-distance functions $\left\{ d_{Y}^{\pi}\right\} $.
Those function does not satisfy the strong projection axiom but they
do satisfy the weak projection axioms, defined below.
\begin{defn}
\label{def:weak-projection-axioms}Let $\by$, $\left\{ {\cal C}Y\right\} $,
$\left\{ \pi_{Y}\right\} $ and $\left\{ d_{Y}\right\} $ be a collection
of projections and its corresponding pseudo-distance function.

We call the set of functions $\left\{ d_{Y}\right\} $ a \textbf{weak
projection} if for some $\theta\geq0$ it satisfy the \textbf{weak}\textbf{\textit{
projection axioms}}:
\begin{itemize}
\item[(P0)] \label{(P0)}$d_{Y}^{\pi}\left(X,X\right)=\diam\pi_{Y}\left(X\right)$
is bounded by $\theta$ for all $X,Y\in\by$;
\item[(P1)]  \label{(P1)}For $X,Y,Z\in\by$ at most one of the numbers
\[
d_{Y}^{\pi}\left(X,Z\right),d_{X}^{\pi}\left(Y,Z\right),d_{Z}^{\pi}\left(X,Y\right)
\]
 is bigger than $\theta$;
\item[(P2)]  \label{(P2)}The set $\left\{ Y\in\by\mid d_{Y}\left(X,Z\right)>\theta\right\} $
is finite;
\end{itemize}

\end{defn}

By definition of the projection sets in HHS it satisfy (P0), from
\cite[1.8]{BHS2} we get (P1) and (P2) is a special case of Masur-Minsky
distance formula for HHS \cite[4.5]{BHS2}. Even though $\mathfrak{L}$,
$\left\{ {\cal C}Y\right\} $, $\left\{ \pi_{Y}\right\} $ and $\left\{ d_{Y}\right\} $
do not satisfy the strong projection axioms, one can coarsely change
$\left\{ \pi_{Y}\right\} $ such that the new projection do satisfy
them.
\begin{thm}[{\cite[Theorem 4.1]{BBFS}}]
 Let $\by$, $\left\{ \left(\mathcal{C}Y,\rho_{Y}\right)\right\} $,
$\left\{ \pi_{Y}\right\} $ and $\left\{ d_{Y}^{\pi}\right\} $ weak
projections with constant $\theta>0$, then projection functions $\left\{ \pi'_{Y}\right\} $
and pseudo-distance $\left\{ d_{Y}^{'}\right\} $ such that
\begin{itemize}
\item $\left\{ d_{Y}\right\} $ satisfy the strong projection axioms for
the constant $11\theta$ ;
\item $\pi'_{Y}\left(X\right)\subseteq N_{\theta}\left(\pi_{Y}\left(X\right)\right)$
;
\item $d_{Y}\left(X,Z\right)=\diam\left(\pi'_{Y}\left(X\right)\cup\pi'_{Y}\left(Z\right)\right)$
(i.e. $d_{Y}$ is obtained from $\pi'$ the same way $d_{Y}^{\pi}$
is obtained from $\pi$) ;
\item If a group $G$ acts on $\by$ preserving the metrics and projection
(i.e. in addition to acting on $\by$ as a set and preserving projections,
we have isometries $\c Y\to\c\left(gY\right)$ for all $Y\in\by$
and $g\in G$) then $G$ also preserve the new projections $\pi'$
(i.e. the construction is $G$-equivariant) ;
\end{itemize}
\end{thm}

This means we can assume that $\mathfrak{L}$, $\left\{ {\cal C}Y\right\} $,
$\left\{ \pi_{Y}\right\} $ and $\left\{ d_{Y}\right\} $ already
satisfy the strong projection (the changing of $\pi$ by $\theta$
does not change any of the argument).

We are now able construct $\projcopx\left(\mathfrak{L}\right)$ and
$\qtms\left(\mathfrak{L}\right)$. Doing so for arbitrary set of pairwise
transverse index won't yield a lot, since if $\Gamma$ is a HHG, in
general there is no action of $\Gamma$ on $\projcopx\left(\mathfrak{L}\right)$
and $\qtms\left(\mathfrak{L}\right)$ since $\mathfrak{L}$ is not
$\Gamma$-invariant, thus we will assume the following property on
our HHG.
\begin{defn}
\label{def:colorable-HHG}A HHG $\Gamma$ is called \textbf{colorable
}if there exists a finite index subgroup $\Gamma'$ such that the
number of $\Gamma'$-orbits in $\mathfrak{S}$ is finite and element
in the same $\Gamma'$-orbits are pair wise transverse.
\end{defn}

Divide $\mathfrak{S}$ into finitely many orbits $\mathfrak{S}=\mathfrak{S}_{1}\cup\ldots\cup\mathfrak{S}_{m}$
as in the definition of colorable. From each $\mathfrak{S}_{i}$ one
can construct a quasi-tree of metric spaces, denoted $\mathcal{C}_{K}^{i}$.
Notice that $\Gamma'$ acts on each $\mathfrak{S}_{i}$ and thus on
each $\c_{K}^{i}$. Since all domains of $\mathcal{X}$ are uniformly
hyperbolic then $\c_{K}^{i}$ is also hyperbolic \cite[4.17]{BBF}.

In addition to the action of $\Gamma$ on each $\c_{K}^{i}$, in order
to apply the techniques in \Secref{Case-of-Product}, we need some
acylindricity assumption.
\begin{defn}
\label{def:stricly-acylidrical-HHG}A colorable HHG $\Gamma$ is called
\textbf{strictly acylindrical} if for every $i=1,\ldots,m$ the action
of $\Gamma'$ on $\qtms^{i}$ is acylindrical.
\end{defn}

As in the case of product diagram, there is no harm in assuming that
$\Gamma'=\Gamma$.

\subsubsection*{Projection to the quasi-trees of metric spaces.}

Assume that $\Gamma$ is a colorable HHG acting on $\mathcal{X}$,
with factors as in the last section. Remember we assume that $\Gamma=\Gamma'$.

In \cite[§5.3]{BBF} the authors define a projection map ${\cal X}\to\qtms^{i}$
in the case of the mapping class group, and in\cite[§3.2]{Ha-Pe}
the authors generalize this construction to a general colorable HHG.

The idea of the projection is that for any $x\in\mathcal{X}$, there
is a coarse point (i.e. a uniformly bounded set) in $\qtms^{i}$,
denoted $\psi_{i}\left(x\right)$, such that for any $U\in\sfrak_{i}$
we have 
\[
d_{\qtms^{i}}\left(\pi_{U}^{\flat}\left(\psi_{i}\left(x\right)\right),\pi_{U}\left(x\right)\right)
\]
 is uniformly bounded independently of $U\in\sfrak^{i}$ and $x$,
where $\pi_{U}^{\flat}$ is the closest point projection $\qtms^{i}\onto\c U$
(where we think of it as a subset of $\qtms^{i}$), and $\pi_{U}:\mathcal{X}\to U$
is the projection from the definition of $\mathcal{X}$ as HHG. Note
that $d_{\qtms^{i}}\left(\pi_{U}^{\flat}\left(\psi_{i}\left(x\right)\right),\pi_{U}\left(x\right)\right)$
is the same as evaluating the distance in $\c U$ itself (\cite[4.2]{BBF})
i.e. $d_{\qtms^{i}}\left(\pi_{U}^{\flat}\left(\psi_{i}\left(x\right)\right),\pi_{U}\left(x\right)\right)=d_{U}\left(\pi_{U}^{\flat}\left(\psi_{i}\left(x\right)\right),\pi_{U}(x)\right)$.
\begin{defn}[Projection to the quasi-trees of metric spaces]
Choose a point in $\mathcal{X}$ marked as $1$(i.e. choose a base
point for $\mathcal{X}$ as the cayley graph of $\Gamma$). Given
$i\in\left\{ 1,\ldots,m\right\} $ there exists (\cite[Lemma 3.5]{Ha-Pe})
a domain $U_{i}\in\sfrak^{i}$ and uniform bound $E$ such that for
all $U\in\sfrak^{i}$ we have $d_{U}\left(1,\rho_{U}^{U_{i}}\right)\leq E$,
i.e. we define $\psi_{i}\left(1\right)=\pi_{U_{i}}\left(1\right)$.
For $x\in\mathcal{X}$ there is $g_{x}\in\Gamma$ ($g_{x}$ is not
unique but the number of such elements is uniformly finite) such that
$d_{\mathcal{X}}\left(g_{x}\cdot1,x\right)\leq D$. We extend $\psi_{i}$
by
\[
\psi_{i}\left(x\right)=g_{x}\pi_{U_{i}}\left(1\right)
\]
We define $\psi:\mathcal{X}\to\prod_{i=1}^{m}\qtms^{i}$ by $\psi\left(x\right)=\left(\psi_{1}\left(x\right),\ldots,\psi_{m}\left(x\right)\right)$.

Note that $\psi$ and $\psi_{i}$ are $\Gamma$-equivariant.
\end{defn}

\begin{lem}[{\cite[Lemma 3.6]{Ha-Pe}}]
$d_{U}\left(\pi_{U}^{\flat}\left(\psi_{i}\left(x\right)\right),x\right)\leq E+D$
for all $U\in\sfrak^{i}$ and all $x\in\mathcal{X}$
\end{lem}

\begin{proof}
For $x=1$ we have that
\[
d_{U}\left(\pi_{U}^{\flat}\left(\psi_{i}\left(1\right)\right),1\right)=d_{U}\left(\pi_{U}^{\flat}\left(\pi_{U_{i}}\left(1\right)\right),1\right)
\]
if $U=U_{i}$ then $\pi_{U}^{\flat}\left(\pi_{U_{i}}\left(1\right)\right)=\pi_{U}\left(1\right)$
so $d_{U}\left(\pi_{U}^{\flat}\left(\pi_{U_{i}}\left(1\right)\right),1\right)=d_{U}\left(1,1\right)\leq E$.
Otherwise $\pi_{U}^{\flat}\left(\pi_{U_{i}}\left(1\right)\right)=\rho_{U}^{U_{i}}$
and by definition of $U_{i}$ we are done.

For other $x\in\mathcal{X}$ we have
\begin{align*}
d_{U}\left(\pi_{U}^{\flat}\left(\psi_{i}\left(x\right)\right),x\right) & \leq d_{U}\left(\pi_{U}^{\flat}\left(g_{x}\psi_{i}\left(x\right)\right),g_{x}\cdot1\right)+D=\\
 & =d_{g_{x}^{-1}U}\left(\pi_{g_{x}^{-1}U}^{\flat}\left(\psi_{i}\left(x\right)\right),1\right)+D\leq E+D
\end{align*}
\end{proof}
The following properties of those projection will be useful in the
rest of the paper
\begin{prop}
\label{prop:qtms-projections}The following holds:
\begin{enumerate}
\item (\cite[3.8]{Ha-Pe}) There is a number $\kappa$ (dependent only of
the HHG constants and $K$) such that $\psi$ is a $\left(\kappa,\kappa\right)$-quasi-isometric
embedding (where the product is endowed with the $l_{1}$ metric).
\item (\cite[3.10]{Ha-Pe}) There are numbers $D,\mu$ (dependent only of
the HHG constants and $K$) such that the image of any $D$-hierarchy
path $\gamma$ under $\psi_{i}$ is an un-parametrized $\mu$-quasi-geodesic
in the hyperbolic space $\qtms^{i}$.
\end{enumerate}
\end{prop}

The following corollary follows directly from the proposition.
\begin{cor}
Each $\psi_{i}$ is coarsely Lipschitz (with constant dependent only
of the HHG constants and $K$).
\end{cor}

\subsection{Higher rank resolution over Colorable HHGs}

Let $\Gamma\acts\mathcal{X}$ be a strictly acylindrical colorable
HHG.

Given a sequence $\left(\varphi_{n}\right)_{n}\in\hom\left(G,\Gamma\right)$
one can look at the sequence of induced action on $\prod_{i=1}^{m}\qtms^{i}$
and continue as in \Subsecref{Higher-rank-Resolution-product}, and
build higher rank resolutions corresponding to action on the product
of the quasi-tree of metric spaces.

In order to mimic the proof of noetherianity in the case of product
diagram, we need an equivalent of \ref{lem:non-divergent-rel-is-non-divergent},
which turns up to be much more difficult. The proof of the next lemma
is due to Sela, it appears inside the proof of \cite[3.3]{Se23},
when Sela shows that if $\out\left(G\right)$ is infinite then the
corresponding higher-rank Makanin-Razborov diagram is non-trivial.
\begin{lem}
\label{lem:non-divergent-rel-is-non-divergent-general-diagram}Given
a sequence $\left(\varphi_{n}\right)\in\hom\left(G,\Gamma\right)$
the following are equivalent:
\begin{enumerate}
\item $\left(\varphi_{n}\right)$ is divergent relative to $\mathcal{X}$
\item $\left(\varphi_{n}\right)$ is divergent relative to the action of
$\Gamma$ on $\mathcal{C}_{K}^{i}$ for some $i$.
\end{enumerate}
\end{lem}

\begin{proof}
Denote by $S=\left\{ s_{1},\ldots,s_{d}\right\} $ a generating set
of $G$, and assume that $1\in S$.

Assume that $\left(\varphi_{n}\right)$ is non-divergent relative
to $\mathcal{X}$. Let $\left(x_{n}\right)_{n}\in\mathcal{X}$ be
such that $\lim_{\omega}\max_{s\in S}d_{\mathcal{X}}\left(\varphi_{n}\left(s\right)x_{n},x_{n}\right)<D$.
Look at the induced action on $\prod_{i=1}^{m}\qtms^{i}$. Since $\psi$
is quasi-isometric embedding it means that $\left(\varphi_{n}\right)_{n}$
is non-divergent relative to $\prod_{i=1}^{m}\qtms^{i}$ (the witnesses
would be the sequence $\left(\psi\left(x_{n}\right)\right)_{n}$ ),
and as in \lemref{non-divergent-rel-is-non-divergent} we get that
$\left(\varphi_{n}\right)_{n}$ is non divergent relative to every
$\qtms^{i}$.

Assume that $\left(\varphi_{n}\right)_{n}$ is non-divergent relative
to all $\qtms^{i}$, but yet divergent relative to $\mathcal{X}$.
Denote the bound of the displacement distances of $\left(\varphi_{n}\right)$
on the $\qtms^{i}$ by $b$.

Our goal is to find $x_{n}\in\mathcal{X}$ such that $\lim_{\omega}\max_{s\in S}d_{\mathcal{X}}\left(\varphi_{n}\left(s\right)x_{n},x_{n}\right)<\infty$ 

\subsubsection*{Step 0}

For each $n$ take some $x_{n}^{0}\in\mathcal{X}$. If $\lim_{\omega}\max_{s\in S}d_{\mathcal{X}}\left(\varphi_{n}\left(s\right)x_{n},x_{n}\right)<\infty$
then we are done, assume that $\lim_{\omega}\max_{s\in S}d_{\mathcal{X}}\left(\varphi_{n}\left(s\right)x_{n},x_{n}\right)=\infty$.

Denote by $A_{n}^{0}\subset\mathcal{X}$ the union of hierarchy paths
$\left[\varphi_{n}\left(s_{i}\right)x_{n}^{0},\varphi_{n}\left(s_{j}\right)x_{n}^{0}\right]$
for all $1\leq i<j\leq d$. For each $i$, we can project $A_{n}^{0}$
to $\qtms^{i}$ using $\psi_{i}$, since $\psi_{i}$ map hierarchy
paths to quasi-geodesics (\Propref{qtms-projections}) and since $\qtms^{i}$
is $\delta$-hyperbolic the image $\psi_{i}\left(A_{n}^{0}\right)$
is $\delta'$ close to a quasi-isometrically embedded finite simplicial
tree, denoted $T_{i,n}^{0}$ ($\delta'$ depends only on $d$ and
$\delta$).

For each $i$ the tree $T_{i.n}^{0}$ is finite and $\lim_{\omega}\diam T_{i.n}^{0}$
might be finite or infinite. We can rearrange the quasi-trees such
that for $1\leq i\leq m'$ we have $\lim_{\omega}\diam T_{i.n}^{0}=\infty$
and for $m'<i\leq m$ we have $\lim_{\omega}\diam T_{i.n}^{0}<\infty$.

From the distance formula, we must have that for at least one index
$\lim_{\omega}\diam T_{i.n}^{0}=\infty$ (since we assumed that $\left(\varphi_{n}\right)$
is divergent relative to $\mathcal{X}$ so $\lim_{\omega}\diam A_{n}^{0}=\infty$),
thus $0<m'$.

\subsubsection*{Step 0.5}

Let's restrict our attention to $i=1$, i.e. to $\qtms^{1}$ and $T_{1,n}^{0}$.
Since $\left(\varphi_{n}\right)_{n}$ is non-divergent relative to
$\qtms^{1}$ then for $\omega$-all values of $n$ there exists $y_{n}\in\qtms^{1}$
such that the displacement distance of the generators on $y_{n}$
is less then $b$. Take $t_{n}^{0}$ to be the closest point projection
of $y_{n}$ on $T_{1,n}^{0}$.

By the choice of $t_{n}^{0}$ as closest point, given $s\in S$, the
paths $\left[y_{n},\varphi_{n}\left(s\right)x_{n}\right]$ and $\left[\varphi_{n}\left(s\right)y_{n},\varphi_{n}\left(s\right)x_{n}\right]$
coarsely (i.e. up to $2\delta$ neighborhood of the paths) contain
$t_{n}^{0}$ and $\varphi_{n}\left(s\right)t_{n}^{0}$ respectively.
Using $\delta$-thin triangle arguments on the triangle $\left[\varphi_{n}\left(s\right)x_{n},y_{n},\varphi_{n}\left(s\right)y_{n}\right]$
one gets that $t_{n}^{0}$ is displaced no more than $b'=b+12\delta$
by $s$, meaning $t_{n}^{0}$ is displaced no more than $b'$ by any
$s\in S$.

Remember that $t_{n}^{0}$ is $\delta'$ close to the projection of
$A_{n}^{0}$ to $\mathcal{\qtms}^{1}$, thus we can coarsely lift
$t_{n}^{0}$ to $A_{n}^{0}$ (in some non-unique way), denoted $x_{n}^{1}$.

\subsubsection*{Step 1}

Now we can repeat the process for $x_{n}^{1}$, to get $A_{n}^{1}$,
$T_{i,n}^{1}$.

Notice that if some $i$ satisfied that $\lim_{\omega}\diam T_{i.n}^{0}<\infty$
then we also have $\lim_{\omega}\diam T_{i.n}^{1}<\infty$. Indeed,
let $i$ be an index such such $\lim_{\omega}\diam T_{i.n}^{0}<\infty$
(i.e. $m'<i$). This is the same as saying that translation of $\psi_{i}\left(x_{n}^{0}\right)$
by the generators has a $\omega$-uniform bound. We took $x_{n}^{1}$
from $A_{n}^{0}$, which coarsely project onto $T_{i,n}^{0}$, thus
the distances between $\psi_{i}\left(x_{n}^{0}\right)$ and $\psi_{i}\left(x_{n}^{1}\right)$
in $\qtms^{i}$ are $\omega$-bounded. Combine the last two sentences,
we get that all the translation of $\psi_{i}\left(x_{n}^{1}\right)$
by the generators has a $\omega$-uniform bound, i.e. $\lim_{\omega}\diam T_{i.n}^{1}<\infty$.

Also, we have that $\lim_{\omega}\diam T_{1,n}^{1}$ is now finite
(since $\psi_{i}\left(x_{n}^{1}\right)$ is $\delta'$ close to $t_{n}^{0}$,
and we took $t_{n}^{0}$ to have bounded translation by the generators.).

This means that when repeating the process with $x_{n}^{1}$ we have
at least one more index such that $\lim_{\omega}\diam T_{i.n}<\infty$.\\

Doing this at most $m$ times we will get series of points $\left(x_{n}^{m}\right)_{n}\in\mathcal{X}$
such that the projections of $A_{n}^{m}$ to each $\qtms^{i}$ are
$\wsurelu$ bounded, meaning that $A_{n}^{m}$ are $\wsurelu$ bounded,
and thus $\left(\varphi_{n}\right)_{n}$ is $\mathcal{X}$-non-divergent.

\end{proof}
\begin{rem}
\lemref{non-divergent-rel-is-non-divergent-general-diagram} and its
proof hold in a more general setting. More precisely, one can replace
$\mathcal{X}$ by any (geodesic) metric space and $\qtms^{i}$ by
any hyperbolic space, as long as $\psi$ and $\psi_{1},\ldots\psi_{m}$
satisfy the conditions in \propref{qtms-projections}.
\end{rem}

\subsection{Proof of Noetherianity}

Now that we have the tools and lemmas from the previous section, one
can replace the term ``factor'' used in \Secref{Case-of-Product}
with the notion of quasi-trees of metric spaces and get equivalent
results to \lemref{base-induction-product}, \propref{step-induction-product},
and \Thmref{product-diagram-noetherian}.

\restate{main-thm}
\begin{proof}
Let $\Gamma'\leq\Gamma$ be the finite-index subgroup in the definition
of colorable HHG. Let $\qtms^{1},\ldots,\qtms^{m}$ be the corresponding
quasi-trees of metric spaces.

We have the action of $\Gamma'$ on $\prod_{i=1}^{m}\qtms^{i}$. Replacing
$X$ with $\prod_{i=1}^{m}\qtms^{i}$ one can apply \propref{step-induction-product}
as is (note that in the proof of this proposition, we haven't used
the assumption that the action is geometric). Thus we reduce to the
case of a sequence $\left(\varphi_{n}\right)\in\hom\left(G,\Gamma'\right)$
which is non-divergent relative to $\prod_{i=1}^{m}\qtms^{i}$. \Lemref{non-divergent-rel-is-non-divergent-general-diagram}
tells us that $\left(\varphi_{n}\right)$ is non-divergent relative
to $\mathcal{X}$. Using \Lemref{base-induction-product} we conclude
that $\left(\varphi_{n}\right)$ $\wsurelu$ factor through the corresponding
limit group (indeed, in \Lemref{base-induction-product} there was
no use of any structure of $\mathcal{X}$, only of the fact that the
action of $\Gamma$ on $\mathcal{X}$ is geometric). Thus $\Gamma'$
is equationally noetherian. Meaning $\Gamma$ is virtually noetherian.
As before, virtually noetherian groups are noetherian themselves (\cite{BMR}),
thus $\Gamma$ is noetherian.
\end{proof}

\bibliography{mybib}{}
\bibliographystyle{alpha}

\end{document}